\documentclass[12pt]{article}
\usepackage{amsfonts,amsmath,graphicx}

\def \beq {\begin{eqnarray}}
\def \eeq {\end{eqnarray}}
\def \beqn {\begin{eqnarray*}}
\def \eeqn {\end{eqnarray*}}

\newcommand{\halmos}{\rule{1ex}{1.4ex}}

\newcounter{for}[section]

\numberwithin{equation}{section}

\newtheorem{itlemma}{Lemma}[section]
\newtheorem{itproposition}[itlemma]{Proposition}
\newtheorem{theorem}[itlemma]{Theorem}
\newtheorem{itcorollary}[itlemma]{Corollary}
\newtheorem{itremark}[itlemma]{Remark}
\newtheorem{itremarks}[itlemma]{Remarks}
\newtheorem{itdefinition}[itlemma]{Definition}
\newtheorem{itexample}[itlemma]{Example}

\newenvironment{fact}{\begin{itfact}\rm}{\end{itfact}}
\newenvironment{claim}{\begin{itclaim}\rm}{\end{itclaim}}
\newenvironment{lemma}{\begin{itlemma}}{\end{itlemma}}
\newenvironment{remark}{\begin{itremark}\rm}{\end{itremark}}
\newenvironment{remarks}{\begin{itremarks} \rm}{\end{itremarks}}
\newenvironment{corollary}{\begin{itcorollary}}{\end{itcorollary}}
\newenvironment{proposition}{\begin{itproposition}}{\end{itproposition}}
\newenvironment{definition}{\begin{itdefinition}\rm}{\end{itdefinition}}
\newenvironment{example}{\begin{itexample}\rm}{\end{itexample}}
\newenvironment{proof}{\noindent {\em Proof}.\ \
}{\hspace*{\fill}$\halmos$\medskip}
\newcommand{\be}[1]{\addtocounter{for}{1} \begin{equation}\label{#1}}
\newcommand{\ee}{\end{equation}}
\newcommand{\bl}[1]{\begin{lemma}\label{#1}}
\newcommand{\br}[1]{\begin{remark}\label{#1}}
\newcommand{\brs}[1]{\begin{remarks}\label{#1}}
\newcommand{\bt}[1]{\begin{theorem}\label{#1}}
\newcommand{\bd}[1]{\begin{definition}\label{#1}}
\newcommand{\bp}[1]{\begin{proposition}\label{#1}}
\newcommand{\bc}[1]{\begin{corollary}\label{#1}}
\newcommand{\bfact}[1]{\begin{fact}\label{#1}}
\newcommand{\bex}[1]{\begin{example}\label{#1}}
\newcommand{\ec}{\end{corollary}}
\newcommand{\efact}{\end{fact}}
\newcommand{\eex}{\end{example}}
\newcommand{\el}{\end{lemma}}
\newcommand{\er}{\end{remark}}
\newcommand{\ers}{\end{remarks}}
\newcommand{\et}{\end{theorem}}
\newcommand{\ed}{\end{definition}}
\newcommand{\ep}{\end{proposition}}
\newcommand{\epr}{\end{proof}}
\newcommand{\bpr}{\begin{proof}}
\newcommand{\bcl}[1]{\begin{claim}\label{#1}}
\newcommand{\ecl}{\end{claim}}

\newcommand{\ecs}{\end{corollary}}
\newcommand{\eers}{\end{exercise}}
\newcommand{\eexs}{\end{example}}
\newcommand{\eems}{\end{example}}
\newcommand{\els}{\end{lemma}}
\newcommand{\eles}{\end{lemmaex}}
\newcommand{\ets}{\end{theorem}}
\newcommand{\eds}{\end{definition}}
\newcommand{\eps}{\end{proposition}}

\newcommand{\bi}{\begin{itemize}}
\newcommand{\ei}{\end{itemize}}
\newcommand{\ben}{\begin{enumerate}}
\newcommand{\een}{\end{enumerate}}

\def\vbar{\mathchoice{\vrule height6.3ptdepth-.5ptwidth.8pt\kern-.8pt}
   {\vrule height6.3ptdepth-.5ptwidth.8pt\kern-.8pt}
   {\vrule height4.1ptdepth-.35ptwidth.6pt\kern-.6pt}
   {\vrule height3.1ptdepth-.25ptwidth.5pt\kern-.5pt}}
\def\fudge{\mathchoice{}{}{\mkern.5mu}{\mkern.8mu}}
\def\bbc#1#2{{\rm \mkern#2mu\vbar\mkern-#2mu#1}}
\def\bbb#1{{\rm I\mkern-3.5mu #1}}
\def\bba#1#2{{\rm #1\mkern-#2mu\fudge #1}}
\def\bb#1{{\count4=`#1 \advance\count4by-64 \ifcase\count4\or\bba A{11.5}\or
   \bbb B\or\bbc C{5}\or\bbb D\or\bbb E\or\bbb F \or\bbc G{5}\or\bbb H\or
   \bbb I\or\bbc J{3}\or\bbb K\or\bbb L \or\bbb M\or\bbb N\or\bbc O{5} \or
   \bbb P\or\bbc Q{5}\or\bbb R\or\bbc S{4.2}\or\bba T{10.5}\or\bbc U{5}\or
   \bba V{12}\or\bba W{16.5}\or\bba X{11}\or\bba Y{11.7}\or\bba Z{7.5}\fi}}

\def \T {{\cal{T}}}
\def \I {{\mathsf{I}}}

\def \A {{\mathcal{A}}}
\def \BB {{\mathcal{B}}}

\def \H {{\cal{H}}}
\def \EE {{\mathcal{E}}}

\def \C {{\cal{C}}}

\def \S {{\cal{S}}}

\def\sqr#1#2{{\vcenter{\vbox{\hrule height .#2pt
                             \hbox{\vrule width .#2pt height#1pt \kern#1pt
                                   \vrule width .#2pt}
                             \hrule height .#2pt}}}}

\def\pmb#1{\setbox0=\hbox{#1}%
   \kern-.025em\copy0\kern-\wd0
   \kern.05em\copy0\kern-\wd0
   \kern-.025em\raise.0433em\box0 }
\def\sqr#1#2{{\vcenter{\vbox{\hrule height.#2pt
     \hbox{\vrule width.#2pt height#1pt \kern#1pt
   \vrule width.#2pt}\hrule height.#2pt}}}}

\def\B{{\mathbb B}}
\def\D{{\mathcal D}}
\def\n{\rho}
\def\N{{\mathbb N}}
\def\Z{{\mathbb Z}}
\def\R{{\mathbb R}}

\def\G{\Gamma}

\def\bs{\backslash}

\def\MV{{\rm MV}}

\def\deg{{\rm deg}}
\def\var{{\rm var}}

\def\reff#1{(\ref{#1})}
\def \ind {\hbox{1\hskip -3pt I}}
\newcommand {\cro}[1] {\left[ {#1} \right]}
\newcommand {\acc}[1] {\left\{ {#1} \right\}}
\newcommand {\pare}[1] {\left( {#1} \right)}

\newcommand {\sous}[1] {\underline{#1}}

\begin{document}
\title{Fluctuations for internal DLA on the Comb.}
\author{
  \renewcommand{\thefootnote}{\arabic{footnote}}
  Amine Asselah \& Houda Rahmani\\
\sl Universit\'e Paris-Est}
\date{}

\footnotetext{
    LAMA, Universit\'e Paris-Est Cr\'eteil,
    61 avenue du g\'en\'eral de Gaulle,
    94010 Cr\'eteil cedex, France\\
    \indent\indent e-mail: amine.asselah@u-pec.fr
    \indent\indent e-mail: houda.rahmani@u-pec.fr
  }
\pagestyle{myheadings}
\markboth
    {Fluctuations for internal DLA on the Comb}
    {Fluctuations for internal DLA}

\maketitle
\begin{abstract}
We study internal diffusion limited aggregation (DLA)
on the two dimensional comb lattice.
The comb lattice is a spanning tree of the euclidean lattice,
and internal DLA is a random growth model,
where simple random walks, starting one at a time
at the origin of the comb, stop when reaching the first 
unoccupied site. An asymptotic shape is suggested
by a lower bound of Huss and Sava \cite{HS-11}.
We bound the fluctuations with respect to this shape.
\end{abstract}

\section{Introduction}\label{sec-intro}
The comb lattice, denoted $\C$, 
is an inhomogeneous spanning tree of $\Z^2$.
Sites $z=(x,y)$ and $z'=(x',y')$ share an edge if either $x=x'$
and $|y-y'|=1$, or if $y=y'=0$ and $|x-x'|=1$; in this case,
we say that $z$ and $z'$ are neighbors. 

Internal DLA on $\C$ is a Markov chain on
the finite subsets of the comb, with initial condition
the empty set, and growing as follows. 
Assume we have obtained a cluster $A$. To build
the cluster with one more site,
launch a simple random walk, on the comb and starting at the origin.
Stop the random walk when it exits $A$, say on site $z$. 
The new cluster is the union of
$A$ and $z$, the first visited site outside $A$. The random
walk with the aggregation rule is called an explorer. We say
that the explorer settles on $z$.

Internal DLA has been first studied on the cubic lattice $\Z^d$
in $d$ dimension. Diaconis and Fulton \cite{DF-01} introduced
it, as well as many variants, with a special emphasize on the
invariance of the cluster with respect to the order in which the
explorers are sent: the so-called {\it abelian property}. Lawler, Bramson
and Griffeath \cite{lawler92} established that the normalized
asymptotic shape is the euclidean sphere in dimension two
or more. 

Blach\`ere and Brofferio \cite{BB-07} obtained
a limiting shape when the graph is a
finitely generated group with exponential growth.
Huss \cite{huss} studied internal DLA for
a large class of random walks on such graphs.
Recently, internal DLA has been considered
on the infinite percolation cluster, and
the asymptotic shape is a euclidean ball intersected with the
infinite cluster: E.Shellef \cite{shellef} obtained
a bound on the inner fluctuations, 
and Duminil-Copin, Lucas, Yadin and Yehudayoff
\cite{Duminil-Lucas} obtained the corresponding 
bound on the outer fluctuations using the inner bound.

\begin{figure}[htpb]
\centering
\includegraphics[width=5cm,height=7cm]{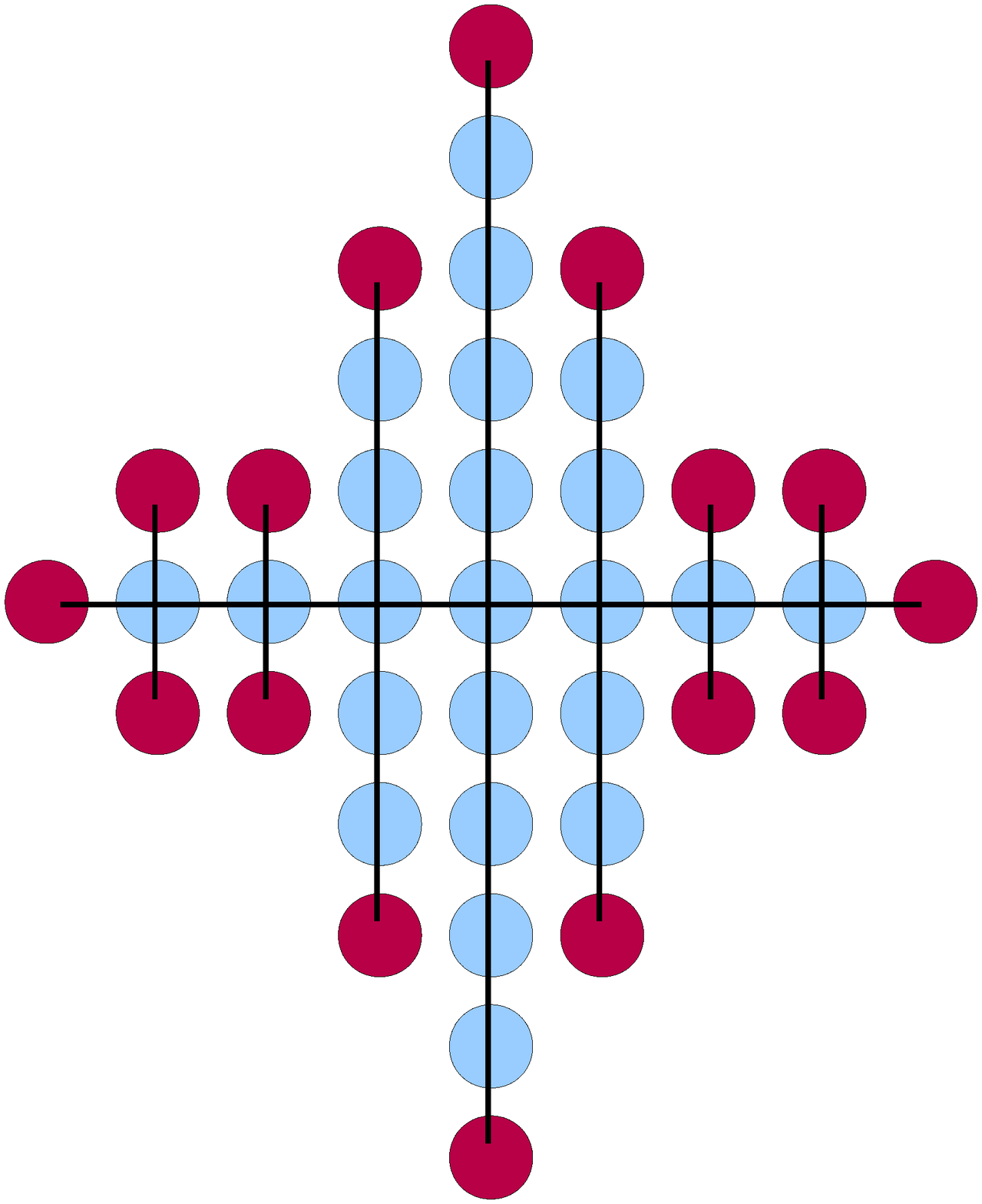}
\caption{$\D(\rho)$ in blue for $\rho=3,5$. Boundary in red.}
\label{fig:comb-cluster}
\end{figure}

It is interesting to study internal DLA on the comb, since it is
inhomogeneous, and distinct from a cubic lattice: 
a simple random walk is recurrent,
however two random walks meet on the average
a finite number of times \cite{KP}.
Also, the $x$ and $y$ axes play
a different role: in a time $n$, a simple random walk on the comb,
reaches a $y$-axis displacement of order $n^{1/2}$, and
an $x$-axis displacement of order $n^{1/4}$.
To discuss results on the comb, let us introduce some notations.
For any real $\rho$, we define 
\be{def-D}
D(\rho)=\acc{(x,y)\in \R^2:\ |x|< \rho,\ 
|y|< \frac{1}{3}\big(\rho-|x|\big)^2},\quad
\text{and}\quad \D(\rho)=D(\rho)\cap \Z^2.
\ee
See Figure~\ref{fig:comb-cluster}.
For an integer $n$, the number of sites in $\D(n)$ is
denoted $d(n)$, and the internal DLA cluster obtained
by sending $d(n)$ explorers is denoted $A(n)$. 

Recently, Huss and Sava~\cite{HS-11} have characterized the
limiting shape for a related model, the {\it divisible sandpile}
introduced in \cite{LP-07}, and shown a lower bound for
the shape of the internal DLA cluster on the comb.
\bt{theo-HS}[Theorem 4.2 of \cite{HS-11}]
For any $\epsilon>0$, with probability 1, we have for
$n$ large enough $\D(n-\epsilon n)\subset A(n)$.
\et
Our main result here is the following improvement.
\bt{theo-amin}
There is a positive constant $a$ such that with probability 1,
and $n$ large enough 
\be{theo-estimates}
\D\big(n-a\sqrt{\log(n)}\big) 
\subset A(n)\subset \D\big(n+a\sqrt{\log(n)}\big).
\ee
\et
\br{rem-main}
This result does not mean that fluctuations are sub-logarithmic, but
rather suggests that they are gaussian. Indeed, a site
$z=(x,y)$ on the boundary of $\D\big(n-a\sqrt{\log(n)}\big)$
is at a distance of order $\frac{2}{3}a\sqrt{\log(n)}(n-|x|)$ from 
the boundary of $\D(n)$, whereas the tooth's length
is of order $\frac{1}{3}(n-|x|)^2$. Thus, the fluctuations
are similar to what would be observed in a collection of $n$
independent segments whose lengths decrease quadratically. 
Diaconis and Fulton in \cite{DF-01} 
used an urn-representation to obtain 
a central limit theorem on $\Z$
for the right-end of the DLA cluster.
\er
Theorem~\ref{theo-amin} follows a classical approach by
Lawler, Bramson and Griffeath~\cite{lawler92}, and
requires a study of the restricted Green's function on $\D$. 
It relies also on a deep connection with another cluster growth model,
{\it the divisible sandpile}, which was discovered by Levine and Peres
\cite{LP-07}. Finally, the limiting shape of
the divisible sandpile cluster was
shown to be $\D(n)$ on the comb in \cite{HS-11}. 

It is interesting to note that exit probabilities from the DLA cluster
are not uniform, as it is the case for the cubic lattice, 
or for discrete groups having exponential growth \cite{BB-07}, or for
the layered  square lattice \cite{KL}. To better appreciate
the following estimate, note that for any $\rho>0$, the size of
the boundary of $\D(\rho)$ is of order $\rho$ 
(see Figure~\ref{fig:comb-cluster}).

\bp{prop-rough} For any real $\rho$, and $z=(x,y)$ in the boundary
of $\D(\rho)$ with $x<\rho$, we have
\be{ineq-green}
\frac{1}{2} \frac{\rho-|x|}{\rho^2+\frac{1}{3}}\le
P_0\big(\text{ The walk exits $\D(\rho)$ at }z \big)
\le \frac{\rho+1-|x|}{\rho^2+\frac{1}{3}}.
\ee
\ep
However, one important property, which holds also for the
cubic lattice, is the following {\it uniform hitting property}.
\bp{prop-uniform}
For each $\rho>0$, there is a stopping time
$\tau^*_\rho$ and two positive constants $\sous \kappa,\bar \kappa$,
independent on $\rho$ such that 
\be{ineq-uniform}
\forall z\in \D(\rho),\qquad\qquad
\frac{\sous\kappa}{|\D(\rho)|}\le P_0\big(S(\tau^*_\rho)=z\big)
\le \frac{\bar\kappa}{|\D(\rho)|}.
\ee
\ep
Proposition~\ref{prop-uniform} is crucial in proving
some large deviation estimates about the cluster,
which in turn, are key in controlling the outer error.

We now turn to some large deviations estimates which
shed more light on the covering mechanism.
Note that a general feature which emerges from all studies
is that during the covering process,
explorers do not leave {\it holes deep inside the bulk}.
Our first lemma deals with the probability an explorer
reaches site $(R,0)$, on the $x$-axis, without leaving
the explored region $V\ni (R,0)$. 
The result, and its proof, are interesting on their own, and
follow closely Lemma 1.6 of \cite{AG2}.
For a subset $\Lambda$ in $\Z^2$, let $H(\Lambda)$ be 
the first time the random walk hits $\Lambda$.
\bl{lem-cross} Let $R$ be a positive integer, and $V$
a subset of $\Z^2$ containing $(0,0)$ and $(R,0)$. 
There are positive constants $a_3$ and $\kappa_3$,
independent of $R$ and $V$, such that
\be{ineq-cross}
P_0\big( H(R,0)<H(V^c)\big)\le \exp\Big(a_3-\kappa_3
\sqrt{ \frac{R^3}{|V|}}\Big).
\ee
\el
In other words, the random walk cannot reach $(R,0)$
inside $V$, unless $V$ contains of the order of $R^3$ sites. 
As a corollary of Lemma~\ref{lem-cross}, we have the
following large deviation upper bound. Recall that $A(n)$ is
the cluster obtained when sending $|\D(n)|$ explorers at the origin,
and that the volume of $\D(n)$ is of order $n^3$.
In other words, Lemma~\ref{lem-cross} quantifies
the probability of making thin tentacles along the $x$-axis.

\bc{cor-cross} There is $\beta,\kappa_2>0$, such that if $R$ and $n$ are
integers with $d(n)<\beta R^3$, then 
\be{ineq-LD}
P\big( \small{(R,0)}\in A(n)\big) \le \exp\big(- \kappa_2 R^2\big).
\ee
\ec

We wish now to explain how to find the asymptotic
shape of the internal DLA cluster built with simple random walks
all started at a distinguished vertex of a graph, say 0.
A fundamental observation of Lawler, Bramson
and Griffeath \cite{lawler92} yields a recipe:
find an increasing family of subsets of the graph,
say $\{\D(\rho),\ \rho\in \R\}$ containing 0, such that the 
discrete {\it mean value property holds} for harmonic functions.
More precisely, $h$ is harmonic if for any vertex $x$ of the graph
\be{def-harmonic}
\sum_{y\text{ neighbor of }x}\big(h(y)-h(x)\big)=0.
\ee
Define now, for any subset $\Lambda$, and any function $h$ on $\Lambda$,
the centered average (recall that walks start from 0)
\be{def-MVP}
\MV(h,\Lambda)=\sum_{z\in \Lambda} \big(h(z)-h(0)\big).
\ee
Finally, we say that {\it the mean value property holds on $\Lambda$}
when for any $h$ harmonic, $\MV(h,\Lambda)$ is of smaller order
than the volume of $\Lambda$ times $h(0)$.
Thus, we look for subsets $\{\D(\rho),\ \rho\in \R\}$,
such that for any $\rho$, 
we can show that {\it the mean value property holds on $\D(\rho)$}.

The observation of \cite{lawler92} behind
the connection between the DLA cluster and the mean value property 
is as follows. Each site of the
DLA cluster is the settlement of exactly one explorer. 
Thus, paint green the explorers' trajectories until settlement, and
add red independent random walks trajectories, starting one on each site
of the cluster. The color-free trajectories, obtained by concatenating
the end-point of a green strand with the red strand starting there, 
are independent random walks starting from 0. 
In short, green explorers glued to red walkers make independent random 
walks all started at 0. Now, if $\D(n)$ is the shape around which
the DLA cluster fluctuates, then the probability a green
explorer exits $\D(n)$ from site $z$ in its boundary is small
if {\it few deep holes are left} as $\D(n)$ gets covered. This
probability is bounded by the difference between the expected number
of random walks starting at 0 and exiting $\D(n)$ from $z$, 
and the expected number of red walkers exiting $\D(n)$ from $z$, 
with one starting on each site of $\D(n)$. 
This difference is $\MV(h_z,\D(n))$, where
$h_z(y)$ is the probability of exiting $\D(n)$ from $z$ when
the initial position of the walk is $y$, and $y\mapsto h_z(y)$ is a 
harmonic function. The smaller is $\MV(h_z,\D(n))$, the better
is the control of the fluctuation of the cluster (see \reff{lawler-3}).

The discrete mean value property holds for spheres on $\Z^d$, 
as shown by Levine and Peres (Theorem 1.3 of \cite{LP-07} and Lemma 6 of
\cite{JLS1}) who used the divisible sandpile for that purpose,
(this property can also be derived from unpublished
estimates of Blach\`ere see the Appendix of \cite{AG1}).
On the comb, a mean value property for the domain $\D(n)$
is essentially contained in the study of Huss and Sava \cite{HS-11}.
This is the starting point of our study.

Let us mention that it is delicate to estimate the shape
of the divisible sandpile cluster. For instance on the comb of $\Z^3$
(where teeth stand on the two dimensional plane), we did not
succeed in identifying the sandpile cluster.

The rest of the paper is organized as follows. We start with
estimates on restricted Green's functions in Section~\ref{sec-notation}. 
Estimates on the Green's function, Propositions~\ref{prop-hitting} 
and \ref{prop-green} are the key technical novelties here.
In Section~\ref{sec-classical}, we recall
the classical approach developed in \cite{lawler92}. The
mean value property is proved in Section~\ref{sec-mean}.
The large deviations estimate Lemma~\ref{lem-cross} and
Corollary~\ref{cor-cross} are proved in Section~\ref{sec-ld}.
Finally, inner and outer errors are respectively
estimated in Sections~\ref{sec-inner}
and~\ref{sec-outer}. Finally, in the Appendix, we prove 
technical properties of the Green's function, most notably
Proposition~\ref{prop-hitting}.

\section{Preliminaries.}\label{sec-notation}
\subsection{Notation}
The comb, denoted $\C$, is a tree rooted at the origin. Any nonzero site
has a unique parent: that is its neighbor which is closer to the origin.
It is convenient to call $\A(z)$ the parent of $z$.

The discrete boundary of $\D(\rho)$, denoted $\partial \D(\rho)$, consists
of the sites of $\Z^2$ not in $\D(\rho)$, 
but at a distance 1 from $\D(\rho)$.
The internal boundary of $\D(\rho)$, denoted $\partial_I \D(\rho)$,
consists of sites of $\D(\rho)$ at a distance 1 from $\partial \D(\rho)$.
The continuous boundary of $D(\rho)$ is denoted $\partial D(\rho)$,
and is the curve $\{(x,y)\in \R^2:\ |x|\le \rho,\ 
|y|= \frac{1}{3}\big(\rho-|x|\big)^2\}$. The euclidean ball
of center 0, and radius $R$ is denoted $B(R)$, and
\[
\B(R)=B(R)\cap \Z^2,\quad
\text{with}\quad B(R)=\acc{(x,y)\in \R^2:\ x^2+y^2<R}.
\]

For a subset $\Lambda$ of $\Z^2$, we denote by $H(\Lambda)$ the time
at which a simple random walk on the comb first hits $\Lambda$,
and we call $\Lambda^+$ 
the intersection of $\Lambda$ with the positive quadrant.
\subsection{On Green's functions.}\label{sec-green}
Henceforth, we consider a simple random walk on the comb $\C$.
We establish many results on harmonic functions on the
domain $\D(\rho)$. 
To ease the reading, their proofs are postponed to the Appendix.

For a subset $\Lambda$ of $\Z^2$,
let $G_\Lambda$ be Green's function restricted to $\Lambda$. In
other words, for $x,y\in \Lambda$, $G_\Lambda(x;y)$ is the expected
number of visits to $y$ before escaping $\Lambda$, when starting on $x$:
\be{def-G}
G_\Lambda(x;y)=E_x\cro{\sum_{n=0}^{\infty} \ind_{n<H(\Lambda^c)} 
\ind_{\{S(n)=y\}}}.
\ee
We first approximate Green's function restricted to $\D(\rho)$.
For a real $x$, $[x]$ denotes its integer part.

\bl{lem-green} Let $\rho$ be any positive real, and $z=(x,y)\in \D(\rho)$.
Define $h$ as
\be{def-h}
h(z)=\frac{2(\rho-|x|)}{\rho^2+\frac{1}{3}}
\Big( \frac{(\rho-|x|)^2}{3}-|y|\Big),
\ee
and $h^+$ as
\be{def-h+}
h^+(z)=\frac{2([\rho]+1-|x|)}{([\rho]+1)^2+\frac{1}{3}+1}
\Big(\frac{([\rho]+1-|x|)^2}{3}+1-|y|\Big).
\ee
Then
\be{green-asymp}
h(z)\le G_{\D(\n)}(0;z)\le h^+(z).
\ee
Moreover, assume that $z'\in \partial \D(\rho)$, $z=\A(z')$.
If $y\not = 0$, then $\n-|x|>1$ and we have
\be{green-boundary}
\frac{(\rho-|X_z|)}{\rho^2+\frac{1}{3}}\le G_{\D(\rho)}(0;z)\le
\frac{2([\rho]+1-|X_z|)}{\rho^2+\frac{1}{3}},
\ee
whereas if $y=0$, then
\be{green-tip}
\frac{1}{4}
\frac{2(\rho+1-|X_z|)}{\rho^2+\frac{1}{3}}\le G_{\D(\rho)}(0;z)\le
\frac{2([\rho]+1-|X_z|)}{\rho^2+\frac{1}{3}}.
\ee

\el
For simplicity, we denote $\EE(\rho)=H(\D^c(\rho))$ the exit time
from $\D(\rho)$. Our second result is our main technical contribution
in estimating hitting probabilities. This in turn allows
us to establish accurate Green's function estimates.

\bp{prop-hitting} For any positive real $\n$, and
any integer $x$, with $|x|\le \n$, there is a constant $\kappa_a>0$
independent of $\n$
\be{lower-hitting}
P_0\big( H(x,0)< \EE(\n)\big)\ge \kappa_a\pare{\frac{\n-|x|}{\n}}^2.
\ee
\ep
We now estimate the probability of exiting $\D(\n)$ from 
$z\in \partial \D(\n)$. This is equivalent to estimating
Green's function $y\mapsto G_{\D(\n)}(y,\A(z))$, since
by a last passage decomposition
\be{last-p}
P_y\big( S(\EE(\n))=z\big)=\frac{1}{\deg(\A(z))} G_{\D(\n)}(y,\A(z)).
\ee
It is convenient for $w\in \Z^2$, to denote
its two coordinates as $X_w$ and $Y_w$. Also,
let $L_\n(w)$ denote the smallest integer
larger or equal than $\frac{1}{3}(\n-X_w)^2$, and
$sg(x)$ is the sign of $x$.

\bp{prop-green}
Assume $z\in\partial \D(\n)$, and $w\in \D(\n)$.\\
(i) When $0\le X_w<X_z$ or $X_w=X_z$ but
$sg(Y_w)\not= sg(Y_z)$, we have
\be{green-1}
P_w\big( S(\EE(\n))=z\big)\le 
\frac{4}{\kappa_a}\frac{L_\n(w)-Y_w}{L_\n(w)}\times
\frac{(\n+1-X_z)}{(\n-X_w)^2}.
\ee
(ii) When $0\le X_z<X_w$, we have
\be{green-2}
P_w\big( S(\EE(\n))=z\big)\le \frac{4}{\kappa_a}\frac{L_\n(w)-Y_w}{L_\n(z)}\times
\frac{(\n-X_w)}{(\n-X_z)^2}.
\ee
(iii) When $X_w<0\le X_z$, there is a constant $\kappa>0$ such that
\be{green-3}
P_w\big( S(\EE(\n))=z\big)\le\kappa \frac{L_\n(w)-Y_w}{L_\n(w)}\times
\frac{(\n-|X_w|)^3(\n+1-X_z)}{\n^5}.
\ee
(iv) when $X_w=X_z$ and $sg(Y_w)= sg(Y_z)$, we have
\be{green-4}
\frac{1}{2} \frac{Y_w}{L_\n(w)}\le P_w\big( S(\EE(\n))=z\big).
\ee
\ep

Finally, we have the following corollary of Proposition~\ref{prop-green}.
\bc{cor-G2}
Assume that $z\in\partial \D(\n)$.
There are constants $\sous\kappa,\bar \kappa>0$ (independent of
$n$ and $z$) such that 
\be{G2-main}
\sous\kappa\ (\n+1-|X_z|)^2\quad
\le \sum_{w\in \D(\n)} P_{w}\big(S(\EE(\n))=z\big)^2
\le\quad \bar \kappa\  (\n+1-|X_z|)^2.
\ee
\ec
\subsection{On a classical approach.}\label{sec-classical}
Denote by $W(\eta, z)$ (resp. $M(\eta,z)$) the number of explorers 
(resp. random walkers) starting on configuration $\eta\in \N^\C$ which hit
$z$. Two special initial configurations play a key role in
internal DLA: we call $d(n)\ind_0$ the configuration
with $d(n)$ explorers at 0; when $\Lambda$ is a subset of $\Z^2$,
we still use $\Lambda$, rather than $\ind_{\Lambda}$, to denote
the configuration with one explorer on each site of $\Lambda$. 
The main observation of \cite{lawler92} yields
the following inequality in law 
\be{main-lawler}
W(d(n)\ind_0,z)+M(A(n),z)\ge M(d(n)\ind_0,z).
\ee
An important feature of \reff{main-lawler} is that $W(d(n)\ind_0,z)$ is
expressed as a difference of two sums of Bernoulli variables.
However, $A(n)$ is unknown, and as such \reff{main-lawler} 
is of little use.
Since we want to show that $A(n)$ is close to a deterministic region
$\D(n)$, we first look for a
region $\I(n)\subset \D(n)$ which is very likely 
covered by the cluster $A(n)$, when $n$ is large. We even require
that $\I(n)$ be covered by explorers not exiting $\I(n)$,
and we call $A_{\I(n)}(n)$ the cluster made by these explorers.
The possibility to discard trajectories exiting $\I(n)$ is
made possible by a key observation of 
Diaconis and Fulton~\cite{DF-01} 
named {\it the abelian property}: the law
of the cluster is independent on the order in which explorers
are launched; this allows to obtain a smaller cluster
if we discard some trajectories.
The key point now is that by definition
\be{bounded-cluster}
A_{\I(n)}(n)\subset \I(n).
\ee
Now, for $z\in \I(n)$, $W_{\I(n)}(\eta,z)$ 
(resp. $M_{\I(n)}(\eta,z)$) denotes the number of explorers
(resp. walkers starting on $\eta$) which hit $z$ before
exiting $\I(n)$. When $z\in \partial\I(n)$, 
$W_{\I(n)}(\eta,z)$ (resp. $M_{\I(n)}(\eta,z)$) still denotes
the number of explorers
(resp. walkers starting on $\eta$) which exit $\I(n)$ from $z$.
The same idea leading to \reff{main-lawler} yields for $z\in \I(n)$
\be{lawler-1}
W_{\I(n)}(d(n)\ind_0,z)+M_{\I(n)}(A_{\I(n)},z)\ge M_{\I(n)}(d(n)\ind_0,z),
\ee
and this inequality becomes an equality when $z\in\partial \I(n)$.
Using \reff{bounded-cluster}, we obtain for $z\in \I(n)$
\be{lawler-2}
W_{\I(n)}(d(n)\ind_0,z)+M_{\I(n)}(\I(n),z)\ge M_{\I(n)}(n\ind_0,z).
\ee
The harmonic function $y\mapsto P_y\big(H(z)\le H(\partial \I(n))\big)$
is denoted $h_z(y)$. Taking the expectation of both sides
of \reff{lawler-2} 
allows a lower bound on the expectation of $W(d(n)\ind_0,z)$
\be{lawler-3}
\begin{split}
E[W(d(n)\ind_0,z)]\ge & E[W_{\I(n)}(d(n)\ind_0,z)]\ge 
\mu(z):=E[M_{\I(n)}(n\ind_0,z)]-
E[M_{\I(n)}(\I(n),z)]\\
=& \big(d(n)-|\I(n)|\big)\times h_z(0)
+MV\big(h_z ,\ \I(n)\big).
\end{split}
\ee
If \reff{lawler-2} were an equality, and using that
$W_{\I(n)}(d(n)\ind_0,z)$ and $M_{\I(n)}(\I(n),z)$ are independent,
we would have the following bound for the
variance of $W_{\I(n)}(d(n)\ind_0,z)$
(we use that the $M$s are sums of Bernoulli)
\be{if-eq}
\var\big(W_{\I(n)}(d(n)\ind_0,z)\big)=\var\big(M_{\I(n)}(n\ind_0,z)\big)
-\var\big(M_{\I(n)}(\I(n),z)\big)\le \mu(z)+\sum_{y\in \I(n)}h^2_z(y).
\ee
Even though \reff{if-eq} is wrong, and that no bound on the
variance is known, \cite{AG1} shows that for a positive constant $\kappa$
\be{ag-1}
P\big(W_{\I(n)}(d(n)\ind_0,z)=0\big)\le \exp\big(-\kappa
\frac{\mu^2(z)}{\mu(z)+\sum_{y\in \I(n)}h^2_z(y)}\big)
\ee
Then, due to the tree structure of the comb $I(n)\not\subset A(n)$
implies that for some $z\in \partial I(n)$, $W_{I(n)}(z)=0$.
We look for a subset $\I(n)$ in $\D(n)$ such that the following
series on the right hand side converges.
\be{lawler-inner}
\begin{split}
\sum_{n\in \N} P\big( \I(n)\not\subset A(n)\big)\le &
\sum_{n\in \N} \sum_{z\in \partial\I(n)} 
P\big( W_{\I(n)}(d(n)\ind_0,z)=0\big)\\
\le & \sum_{n\in \N} \sum_{z\in \partial\I(n)} 
\exp\big(-\kappa
\frac{\mu^2(z)}{\mu(z)+\sum_{y\in \I(n)}h^2_z(y)}\big).
\end{split}
\ee
Using Borel-Cantelli, \reff{lawler-inner} implies that almost surely,
for $n$ large enough, $\I(n)\subset A(n)$. Since
the lower bound of the asymptotic shape theorem proved by
Huss and Sava \cite{HS-11}, implies that for $n$ large enough
$\I(n)\subset A(n)$, we would conclude the proof of Theorem
\ref{theo-amin}.

This approach can be implemented if we can estimate $\mu(z)$, and
the sum of $y\mapsto h^2_z(y)$ over $\I(n)$. 
Note that $\mu(z)$ should be of
order $(|\D(n)|-|\I(n)|) h_z(0)$ provided we can show that
\be{to-show}
\forall z\in \I(n)\qquad \MV(h_z,\I(n))\ll 
\big(|\D(n)|-|\I(n)|\big) h_z(0).
\ee
\paragraph{The divisible sandpile.}
Levine and Peres \cite{LP-07} have introduced a model,
the divisible sandpile, whose cluster is a good candidate for $\D(n)$.
In this model, we start with a mass $n$ at the origin of our graph,
and topple sand along some sequence of sites. We topple the sand
at a site if its mass is above 1, and we transfer the total mass minus
1 equally to each nearest neighbor. The toppling sequence
is arbitrary provided it covers each site of the graph infinitely
often. We call $z\mapsto w_n(z)$ the final sand distribution, and
we call $z\mapsto u_n(z)$ the odometer function: that is the amount
of sand emitted from each site. The {\it sandpile cluster} is
$\S_n=\{z:\ u_n(z)>0\}$. The key observation is that for
any harmonic function on $\S_n$,
\be{lawler-5}
\sum_{z\in \Z^2} w_n(z)\big(h(z)-h(0)\big)=0.
\ee
When the graph is the comb $\C$, Huss and Sava obtain in \cite{HS-11}
the following result.
\bt{theo-huss}[Theorem 3.5 of \cite{HS-11}]
There is a positive constant $R_{HS}$, such that
for $n$ large enough
\be{bound-huss}
\D(n-R_{HS})\subset\S_n\subset \D(n+R_{HS}).
\ee
\et
This result 
Theorem~\ref{theo-huss} of \cite{HS-11}
is not precise enough for our purpose, but
the arguments in \cite{HS-11} yield easily the following stronger result.
When $A,B$ are subsets of $\Z^2$, it is handful to use the notation $A+B$
for Minkowski addition $\{z=x+y:\ x\in A,\ y\in B\}$, and
$A-B=\{z=x-y:\ x\in A,\ y\in B\}$. 
\bl{lem-cluster} There is a constant $R_{HS}>0$, such that for $n$ 
large enough
\be{bound-amin}
\D(n)-\B(R_{HS})\subset\S_n\subset \D(n)+\B(R_{HS}).
\ee
\el
To prove the lemma, it is enough to check that on
the (continuous) boundary of $\D(n)$, the obstacle function $\gamma_n$
is bounded by a constant, independent of $n$. By the symmetry
of $\D(n)$, it is enough to consider
$x,y$ both positive satisfying $x\le n$ and
$y=\frac{1}{3}(n-x)^2$. We recall 
Huss and Sava's expression of $\gamma_n$ with our normalizing of $\D(n)$
(that is if $n'$ is their $n$, then $n^3=9n'/4$):
\be{huss-1}
\gamma_n(x,y)=\frac{1}{2}\Big( y- \frac{1}{2}\big(
\frac{2}{3}x^2-tx+\frac{9}{24}t^2+\frac{1}{6}\big)\Big)^2,
\ee
with (using the value for $T(n)$ after (3.12) of \cite{HS-11}
with our $n$)
\be{huss-2}
t=T-\frac{20}{27}\frac{1}{T},\quad\text{and}\quad
T=\frac{4}{3}n +O(\frac{1}{n^5}).
\ee
Note that
\be{huss-3}
t=\frac{4}{3}n -\frac{5}{9}\frac{1}{n}+O(\frac{1}{n^5}),\quad
\text{and}\quad
t^2= \big(\frac{4}{3}\big)^2n^2-\frac{40}{27}+O(\frac{1}{n^2}).
\ee
A simple computation yields for $0\le x\le n$, and $y=\frac{1}{3}(n-x)^2$
\be{huss-4}
\gamma_n(x,y)=\frac{1}{2}\Big( -\frac{5}{18}\frac{x}{n}+\frac{7}{36}+
O(\frac{x}{n^2})\Big)^2.
\ee
Thus, for a constant $K$ independent of $n$,
\be{huss-5}
\sup_{z\in\partial_c D(n)} \gamma_n(z)\le K.
\ee
Now, the obstacle function is a upper bound for the odometer
$u_n$ which decays by one unit as we move along a tooth (or along the
$x$-axis), away from the origin. Thus, there is $R_{HS}$ such that
$\S_n\subset \D(n)+\B(R_{HS})$.
Note that on $\partial(\D(n)+\B(R_{HS}))$ the
odometer vanishes, whereas $\gamma_n$ is bounded uniformly in $n$,
say by $\tilde K$. 
Since $u_n-\gamma_n$ is superharmonic, it satisfies the minimum
principle, and satisfies in $\D(n)+\B(R_{HS})$ that 
$u_n\ge \gamma_n-\tilde K$. Since $\gamma_n$ increases quadratically
as we move toward the origin, this implies the lower estimate
$\D(n)-\B(R_{HS})\subset\S_n$ for some 
constant $R_{HS}$ independent of $n$.

\subsection{On the mean-value property in $\D(\rho)$}\label{sec-mean}
Our main result in this section is the following mean value
approximation, which relies on Lemma~\ref{lem-cluster}, where
the constant $R_{HS}$ appears. We consider $z\in \partial
\D(\rho)$, and for $y\in \D(\rho)$ we set $h_z(y)=P_y(S(\EE(\rho))=z)$.

\bl{lem-MV} 
There is a constant $C_{\MV}>0$, such that for any $\rho>0$
and any $z\in \partial\D(\rho)$
\be{ineq-MV}
\big|\MV\big(h_z;\D(\rho)\big)\big|\le\ C_{\MV}R_{HS}^2.
\ee
\el
\br{rem-MV}
For the outer fluctuation, we need a related and simpler result,
that we present now. We consider $\rho'<\rho-R_{HS}$, and
have that for some positive constant $C_{\MV}$
\be{MV-outer}
\big|\MV\big(h_z;\D(\rho')\big)\big|\le\ C_{\MV}R_{HS}.
\ee
We explain after the proof of Lemma~\ref{lem-MV} how to obtain this
simpler statement.
\er
\bpr
First, we extend $h_z:\D(\rho)\to[0,1]$
into a harmonic function on the smallest sandpile cluster, say $\S$,
containing $\D(\rho)$.
By Lemma~\ref{lem-cluster}, it is enough to extend it to
$\D(\rho)+\B(R_{HS})$ with the constant $R_{HS}$ appearing there.

We set $\tilde h\equiv h_z$ on $\D(\rho)\cup\partial\D(\rho)$.
Let $w\in \partial\D(\rho)$ with $|X_w|\le \rho-1$. This implies that
\[
h_z(w)=0,\quad\text{and}\quad h_z(\A(w))>0.
\]
Since $w$ is not on the $x$-axis, there is a unique site $w'$ such
that $\A^k(w')=w$, we denote for simplicity $w'=\A^{-k}(w)$.
Since teeth are one-dimensional, harmonicity of $\tilde h$
imposes that for any positive integer $k$
\[
\tilde h(\A^{-k}(w))-\tilde h(w)=k\big(
h_z(w)-h_z(\A(w))\big)=-kh_z(\A(w)),
\]
so that if $w\in \partial\D(\rho)$ but $w$ not on the $x$-axis,
\be{ext-G}
\tilde h(\A^{-k}(w))=-k h_z(\A(w)).
\ee
On the $x$-axis, we choose the following extension
\be{ext-G2}
\tilde h\big([\rho]+k,0\big)=-(k-1) h\big([\rho],0\big)
\ee
Now, 
and if $([\rho]-1,1)\not\in \D(\rho)$ we set 
for $l\in \Z$, we set
\[
\tilde h\big([\rho]-1,l\big)=-(l-1)h_z([\rho]-1,0).
\]
Finally, we note that $([\rho],1)\not\in \D(\rho)$, and we
extend $\tilde h$ by linearity on each tooth rooted 
on $\{([\rho],k),\ k\in \N\}$, so that for integers
$k\ge 0$, and $l\in \Z$, we have
\be{ext-G3}
\tilde h\big([\rho]+k,l\big)= (l-1)(k-1) h\big([\rho],0\big),
\quad\text{and}\quad
\tilde h\big(-[\rho]-k,l\big)= (l-1)(k-1) h\big(-[\rho],0\big).
\ee
Using a result of Levine and Peres (Theorem 1.3 of \cite{LP-07}), 
and Lemma 2.6 of \cite{HS-11},
there exists a function $y\mapsto \omega(y)$ with value in $[0,1]$, which 
vanishes on $\D(\rho)+\B(R_{HS})$, which equals 
1 on $\D(\rho)-\B(R_{HS})$, and which satisfies
\be{MV-Levine}
\sum_{y\in \Z^2}\omega(y)\big( \tilde h(y)-\tilde h(0)\big)=0.
\ee
Thus, if we denote $\partial_R\D(\rho)$
the shell $\D(\rho)+\B(R_{HS})\bs(\D(\rho)-\B(R_{HS}))$, we have
\be{MV-1}
\begin{split}
\big|\sum_{y\in \D(\rho)}\tilde h(y)-\tilde h(0)\big|=&
\big| \sum_{\D(\rho)\bs(\D(\rho)-\B(R_{HS})}
(1-\omega(y))(\tilde h(y)-\tilde h(0))\\
&\qquad -\sum_{\D(\rho)+\B(R_{HS})\bs\D(\rho)}
\omega(y)(\tilde h(y)-\tilde h(0))\big|
\le\sum_{y\in \partial_R\D(\rho)} \big| \tilde h(y)-\tilde h(0)\big|.
\end{split}
\ee
This implies that for some positive constant $C$
\be{MV-2}
\begin{split}
\big| MV\big(h_z(\cdot),\D(\rho)\big)\big|= &
\big| \sum_{y\in \D(\rho)} \big( \tilde h(y)-\tilde h(0)\big)\big|
\le \!\!\!\sum_{y\in \partial_R\D(\rho)}\!\! | 
\tilde h(y)-\tilde h(0)|\\
\le & \big| \partial_R\D(\rho)\big|\times h_z(0)+
C R_{HS}^2 \sum_{y\in \partial \D(\rho)} h_z\big(\A(y)\big)\\
\le & C R_{HS}^2
\Big(\frac{\rho(\rho-X_z)}{\rho^2}+
\sum_{y\in \partial \D(\rho)} h\big(\A(y)\big)\Big).
\end{split}
\ee
The following bound implies \reff{ineq-MV}. It is
a consequence of Proposition~\ref{prop-green}, after
we decompose $\partial \D(\rho)$ into four regions to be dealt with
estimates \reff{green-3}, \reff{green-1}, \reff{green-4}
and \reff{green-2}. Thus, there is a positive constant $K$ such that
\be{MV-3}
\sum_{y\in \partial \D(\rho)} h_z\big(\A(y)\big)\le 
\sum_{k=1}^{[\rho]} \frac{k^3(\rho+1-X_z)}{k^2\rho^5}+
\sum_{k\ge \rho+1-X_z}\!\!\frac{(\rho+1-X_z)}{k^4}+1+
\sum_{k=1}^{[\rho]+1-X_z}\!\!\frac{k}{(\rho+1-X_z)^4} \le K.
\ee
This concludes the proof of Lemma~\ref{lem-MV}. Finally, we wish
to explain Remark~\ref{rem-MV}. First, $h_z$ is harmonic
on the smallest sandpile cluster containing $\D(\rho')$,
so there is no need to extend it as in the previous proof.
The estimates \reff{MV-2} yields here
\be{MV-5}
\big| MV\big(h_z(\cdot),\D(\rho')\big)\big|
\le \sum_{y\in \partial_R\D(\rho')}\!\! | h_z(y)-h_z(0)|
\le \big| \partial_R\D(\rho')\big|\times h_z(0)+
\sum_{y\in \partial_R\D(\rho')} h_z(y).
\ee
It is now enough to note that on each tooth
intersecting $\partial_R\D(\rho')$ there are
at most $2R_{HS}$ sites, and that the estimates
for $h_z(w)$ in Proposition~\ref{prop-green} are 
worse when $Y_w=0$, and this yields the bound
\be{MV-6}
\sum_{w\in \partial_R\D(\rho')} h_z(w)\le 2R_{HS} 
\sum_{|x|\le \rho'+R_{HS}} h_z(x,0)\le C_{MV} R_{HS}.
\ee
\epr

\section{Large Deviations.}\label{sec-ld}
Our aim in this section is to prove Proposition~\ref{prop-uniform},
Lemma~\ref{lem-cross} and Corollary~\ref{cor-cross}.
\subsection{On the uniform hitting property.}\label{sec-uniform}
For each $\rho>0$, we build here a stopping time
$\tau^*_\rho$ which satisfies \reff{ineq-uniform} of
Proposition~\ref{prop-uniform}.
The time $\tau^*_\rho$ is called a {\it flashing time}.

We set
\be{def-grho}
g_\rho(r)=\frac{3r^2}{\rho^3}\quad\text{for }r\in [0,\rho], \quad
\text{and for }r>\rho\quad g_\rho(r)=0.
\ee
The algorithm which defines $\tau^*_\rho$ is as follows.
\begin{itemize}
\item Draw $R$ according to $g_\rho$.
\item If $R<\frac{1}{2}$, then $\tau^*_\rho=0$, and the walk flashes
on its initial position, the origin.
\item If $R> \frac{1}{2}$, then 
$\tau^*_\rho=\inf\{ t>0: S(t)\not\in \D(R)\}$.
\end{itemize}
We need to estimate $P_0(S(\tau^*_\rho)=z)$ for $z\in \D(\rho)$.
We have
\[
P_0(S(\tau^*_\rho)=z)=P(R<\frac{1}{2})\ind_{z=0}+
\int_{1/2}^\rho P_0\big( S(\EE(r))=z\big) g_\rho(r) dr\times
\ind_{z\not = 0}.
\]
First, $P(S(\tau^*_\rho)=0)=P(R<\frac{1}{2})= 1/(2\rho)^3$.
Assume henceforth that $z\not=0$, and note
that $S(\EE(r))=z$ is possible only if $z\in \partial\D(r)$.
Thus, we define $R(z)<\bar R(z)$ such that
\be{def-R}
z\in \partial\D(r)\Longleftrightarrow R(z)<r\le \bar R(z).
\ee
In other words, we define
\be{def-Rbis}
\frac{1}{3}\big(\bar R(z)-|X_z|\big)^2=|Y_z|,\quad\text{and}\quad
R(z)=\bar R(\A(z)).
\ee
We need to estimate $P_0\big( S(\EE(r))=z\big)$ for $R(z)\le r< \bar R(z)$.
Upper and lower bounds are obtain for Green's function 
in Lemma~\ref{lem-green}, and hold for the exit distribution
by the last passage decomposition \reff{last-p}.
Now, the upper and the lower bound for $P_0\big( S(\tau^*_\rho)=z\big)$
are done in a similar way, and we write in details only the upper
bound. Also, because of the symmetry of $\D(\rho)$, we can assume that
$X_z\ge 0$ and $Y_z\ge 0$. We treat three cases: (i) when $z$ is
a nearest neighbor of the origin, (ii) when $\A(z)\not= 0$ and
$z$ is not on the $x$-axis, and (iii) when $\A(z)\not= 0$ and $Y_z=0$.
\paragraph{Case (i): $\A(z)=0$.} Then, $\bar R(z)\le 2$ and $R(z)=0$.
We have
\be{case-i}
\begin{split}
P_0(S(\tau^*_\rho)=z)\le &\int_{1/2}^{\bar R(z)} \frac{1}{\deg(z)}
\frac{2(r+1-X_z)}{r^2}\frac{3r^2}{\rho} dr\\
\le & \frac{3}{2\rho^3}\big(\bar R(z)-\frac{1}{2}\big)
\big(\bar R(z)+1-X_z\big)\le \frac{27}{4\rho^3}.
\end{split}
\ee
\paragraph{Case (ii): $\A(z)\not=0$ and $Y_z\not =0$.}
Note that $\A(z)=(X_z,Y_z-1)$, $R(z)\ge 1$, and
\[
\bar R(z)-R(z)= \sqrt{3Y_z}-\sqrt{3(Y_z-1)}\le 
\frac{\sqrt 3}{\sqrt{Y_z}}.
\]
Then
\be{case-ii}
\begin{split}
P_0(S(\tau^*_\rho)=z)\le &\int_{R(z)}^{\bar R(z)} \frac{1}{\deg(z)}
\frac{2(r+1-X_z)}{r^2}\frac{3r^2}{\rho} dr\\
\le & \frac{3}{2\rho^3}\big(\bar R(z)+1-X_z\big)
\big(\bar R(z)-R(z)\big)\le \frac{3}{2\rho^3} \sqrt{3Y_z}
\frac{\sqrt 3}{\sqrt{Y_z}}\le \frac{9}{2\rho^3}.
\end{split}
\ee
\paragraph{Case (ii): $\A(z)\not=0$ and $Y_z=0$.} Then 
$\A(z)=(X_z-1,0)$, $\bar R(z)=X_z$ and $R(z)=X_z-1\ge 1$. We have
\be{case-iii}
P_0\big(S(\tau^*_\rho)=z\big)\le\int_{X_z-1}^{X_z} \frac{1}{\deg(z)}
\frac{2(r+1-X_z)}{r^2}\frac{3r^2}{\rho} dr
\le \frac{3}{4\rho^3}.
\ee
We omit the similar estimates yielding the lower bound
of \reff{ineq-uniform}. 

\subsection{Proof of Lemma~\ref{lem-cross}.}
We consider here that the explored region is $V$, and estimate
the probability an explorer reaches $(R,0)$. To obtain \reff{ineq-cross} 
we can assume that the ratio $|V|/R^3$ is as small as we wish.
Also, we can restrict to $|V|\ge R$, since
an explorer reaching $(R,0)$ has to visit all sites 
of $\{(x,0),\ 0\le x\le R\}$.

The proof makes use of the concept of {\it flashing explorer},
which was introduced in \cite{AG1}, and follows the arguments
of the proof of Lemma 1.6 of \cite{AG2}. A flashing explorer is
a random walk which settles only if at some times, {\it the flashing times},
it is not on the explored region $V$.

If an explorer reaches $(R,0)$ (without escaping $V$), then a
flashing explorer following the same trajectory
would reach $(R,0)$ as well. Since Lemma~\ref{lem-cross} requires
a bound from below on the crossing probability, it is enough
to obtain an estimate for the flashing explorer.

\begin{figure}[htpb]
\centering
\includegraphics[width=9cm,height=6cm]{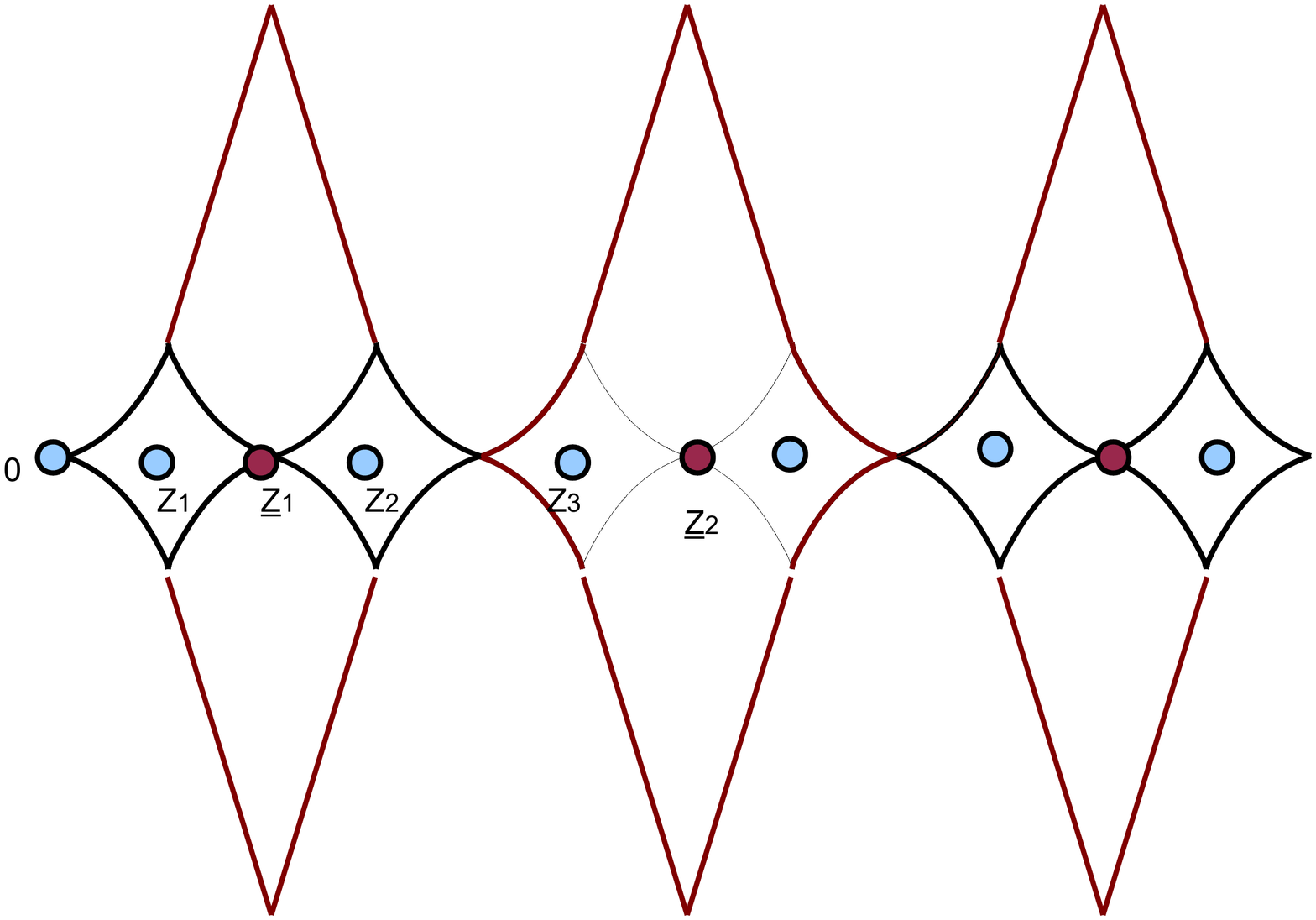}
\caption{Two scales $h$ and $\sous h=2h$, and corresponding centers.}
\label{fig:comb-scale}
\end{figure}

We now define the flashing explorer associated with
{\it the scale} $h$. Let $h$ be a positive real smaller than $R/2$, and
write $M_h$ for the integer part of $R/(2h)$. 
We form $M_h$ disjoint domains
by translating $\D(h)$ so that they cover $[0,2hM_h]$, see
Figure~\ref{fig:comb-scale}, and we call
$Z_1,\dots,Z_{M_h}$ their centers. The flashing explorer associated
with scale $h$ is as follows.
\begin{itemize}
\item It performs a simple random walk on the comb, starting at 0.
\item The first time the walk reaches $Z_i$, it draws one variable
$R_i$ according to $g_h$ and a flashing time $\tau^*_i$ is
constructed as in the previous section but around $\D(Z_i,h)$.
\item It settles the first time $H(Z_i)+\tau^*_i$
that $S(H(Z_i)+\tau^*_i)\not\in V$, for $i=1,\dots,M_h$.
\end{itemize}
We say that the domain $\D(Z_i,h)$ is
{\it well-covered} when $|\D(Z_i,h)\cap V|>\beta |\D(Z_i,h)|$, 
for a positive $\beta<1$ to be chosen later.
We call $\G_h$ the set of well-covered domains: 
\be{def-Gh}
\G_h=\acc{i\in [1,M_h]:\ |\D(Z_i,h)\cap V|>\beta |\D(Z_i,h)|}.
\ee
The reason we use a flashing explorer is that 
the probability that it settles in a {\it not well-covered} domain
is easy to estimate.
Indeed, by the {\it uniform hitting property},
it visits the domain $\D(Z_i,h)$ almost uniformly, and the probability
it flashes on a site of $V$ is less than $\kappa \beta$, for
some positive constant $\kappa<1$ (independent of $\beta$ and $h$).
We now choose $\beta$ by requiring that when $h=R/2$ (and $M_h=1$ and
$Z_1=(h,0)$), then $\G_h=\emptyset$.

Thus, 
\be{cross-7}
P_0\big(\text{ The flashing explorer reaches }(R,0)\big)\le
\prod_{i\not\in \G_h}P_{Z_i}\big(S(H(Z_i)+\tau^*_i)\in V\big)\le
(\kappa \beta)^{M_h-|\G_h|}.
\ee
Now, $\beta |\D(h)| |\G_h|\le |V|$, and for some $\alpha>0$
$|\D(h)|\ge \alpha h^3$.
Replacing $|\G_h|$ in \reff{cross-7} by an upper bound yields
\be{cross-8}
\begin{split}
P_0\big(\text{ The flashing explorer reaches }(R,0)\big)\le&
\exp\Big( -\log(\frac{1}{\kappa \beta})\big( M_h-\frac{|V|}
{\beta|\D(h)|} \big)\Big)\\
\le & \exp\Big( -\log(\frac{1}{\kappa \beta})
\big(\frac{R}{2h}-1 -\frac{|V|}{\beta \alpha h^3} \big)\Big)
\end{split}
\ee
We optimize now in $h$. 
The maximum of $R/(2h)-|V|/(\beta\alpha  h^3)$ is reached
for $h^*=\sqrt{6/\beta\alpha }\sqrt{|V|/R}$, and this choice completes
the proof of \reff{ineq-cross}.

\subsection{Proof of Corollary~\ref{cor-cross}.}
The proof relies on Lemma~\ref{lem-cross}, and
follows closely the arguments of Lemma 1.5 of
\cite{AG2} with $d=3$, since $|\D(R)|$ is of order $R^3$. 

The strategy of the proof is to built optimal
disjoint random domains $\D(Z_0,h_0)$... 
$\D(Z_L,h_L)$ inside $\D(0,R)$ in such a that if
$N_i$ is the number of settled explorers in $\D(Z_i,h_i)$, 
we have that $(R,0)\in A(n)$ implies that for each $j$,
$N_{j+1}+\dots+N_{L}$
explorers have crossed $\D(Z_0,h_0)\cup\dots\cup\D(Z_j,h_j)$.
The randomness comes from $A(n)$.

We choose $h_0=R/4>1$, and $Z_0=(h_0,0)$. 
We choose a positive (large) constant $\gamma$ from $a_3$ and
$\kappa_3$ of Lemma~\ref{lem-cross}:
\[
\gamma=\max\big(1,(\frac{2a_3}{\kappa_3})^2\big).
\]
The choice of $\gamma$ will be clear later.
We choose now $\beta$ such that $\gamma d(n)<\gamma \beta R^3<|\D(h_0)|$.

We now build by induction neighboring domains $\D(Z_i,h_i)$ for 
$i=1,\dots,L$ such that
\[
\sum_{i=0}^L \D(Z_i,h_i)\supset\acc{(0,k):\ 0\le k\le R,\ 
k\in \N}.
\]
Assume we have chosen $h_{i-1}\ge 1$ and $2(h_0+\dots+h_{i-1})<R$. 
We choose $h_i$ such that
\be{old-3}
|\D(h_i)|=\gamma \big|  \D(Z_{i-1},h_{i-1})\cap A(n)\big|,
\quad\text{and}\quad
Z_i=\big(2\sum_{j=0}^{i-1}h_j+h_i,0\big).
\ee
Note that since $h_{i-1}\ge 1$, we have $N_{i-1}=|
 \D(Z_{i-1},h_{i-1})\cap A(n)|\ge 3$, and $h_i\ge 1$.
Clearly, the induction stops with $L$ domains, and $L\le R/2$.

For any choice of integers $l,n_0,\dots,n_l$, the event
$\{L=l,N_0=n_0,\dots,N_l=n_l\}$ implies that $n_1+\dots+n_l$
explorers have crossed $\D(Z_{0},h_{0})$ with an explored region
made up of $n_0$ settled explorers, and $n_2+\dots+n_l$ 
explorers have crossed $\D(Z_{1},h_{1})$ with an explored region
made up of $n_1$ settled explorers, and so on and so forth.
We use now \reff{ineq-cross} of Lemma~\ref{lem-cross} to obtain
\be{old-2}
\begin{split}
P\big(\text{ an explorer reaches }&(R,0)\big)\le
\sum_{l,n_0,\dots,n_l} \prod_{i>0}
P\big(\sum_{k=i}^l n_k\text{ explorers cross }\D(Z_{i-1},h_{i-1})\big)\\
\le (R^3)^R\sup_{l,n_1,\dots,n_l}&
\exp\Big( a_3\sum_{i=1}^l in_i-\kappa_3\sum_{i=1}^l
n_i\big(\sqrt{\frac{h^3_0}{n_{0}}}+\dots+
\sqrt{\frac{h^3_{i-1}}{n_{i-1}}}\big)\Big).
\end{split}
\ee
By the arithmetic-geometric inequality, for $1\le i \leq l$
(and using $h_i \leq h_0$)
\be{old-4}
\begin{split}
\frac{1}{i}\big(
\sqrt{\frac{h_0^3}{n_0}}+\dots+\sqrt{\frac{h_{i-1}^3}{n_{i-1}}}\big)
\ge&\big(\frac{h_0^3}{n_0}\times\dots\times \frac{h_{i-1}^3}{n_{i-1}}
\big)^{1/2i}\\
=&\pare{\frac{h_0^3}{n_{i-1}}\gamma^{i-1}}^{1/2i}
= \pare{\frac{h_0^3}{h_i^3}\gamma^{i}}^{1/2i}
\ge \frac{2a_3}{\kappa_3}.
\end{split}
\ee
Thus, from \reff{old-2} and \reff{old-4}, we have
\be{old-5}
\begin{split}
P\pare{ (R,0)\in A(n)} & \le R (\beta R^3)^{R+1}
\max_{
  \stackrel{\scriptstyle l\le R, n_0,n_1,\dots,n_l\leq \beta R^3}
  {\scriptstyle \forall i, n_i \geq \lfloor h_i\rfloor}
}
\exp\pare{-a_3 \sum_{i=1}^{l} i\ n_{i}}\\
&\le R (\beta R^3)^{R+1}
\max_{
  \stackrel{\scriptstyle l\le R, n_0,n_1,\dots,n_l\leq \beta R^3}
  {\scriptstyle \forall i, n_i \geq \lfloor h_i\rfloor}
}
\exp\pare{-\frac{a_3}{\gamma} \sum_{i=1}^{l-1} i h^3_{i+1}}
\end{split}
\ee
Since $h_1\le R/4$, we have $h_2+\dots+h_l\ge R/4$.
By H\"older's inequality, and for a constant $c_3$, we have
\be{holder}
\sum_{i=1}^{l-1} i h^3_{i+1}\ge 
\frac{\pare{\sum_{i=1}^{l-1} h_{i+1}}^3}{
\pare{\sum_{i=1}^{l-1}\frac{1}{\sqrt i}}^{2}}
\ge c_3 \frac{R^3}{l}\ge c_3 R^2.
\ee
This concludes the proof.

\section{Inner Fluctuations.}\label{sec-inner}
Our main result here is an inner estimate for the aggregate.
\bp{prop-inner}
There is a positive constants $\kappa_{in}$ (independent of $n$) such that
for any positive $a$ and integer $n$ large enough
\be{ineq-inner}
P\big(\D(n-a\sqrt{\log(n)})\not\subset A(n)\big)\le
d(n) \exp\big( -\kappa_{in} a^2\log(n)\big).
\ee
\ep
\br{rem-inner}
Since $d(n)$ is of order $n^3$, \reff{ineq-inner} implies
the inner estimate of \ref{theo-estimates} of Theorem~\ref{theo-amin}.
Our proof below establishes actually a stronger result
than \reff{ineq-inner}. Indeed, we only count explorers which
remain in the domain $\D(n)$. This remark is used in the outer
error bound.
\er
\bpr
The constant $A>1$ will be chosen later. For any $\alpha>0$, 
we set $a=A\sqrt{\alpha}$ and $L=\sqrt{\alpha\log(n)}$. 
Inequality \reff{ineq-inner} follows if we show that for
$z\in \partial\D(n-L)$ with $|X_z|\le n-AL$, we have 
\be{inner-one}
P\big(W_{\D(n-L)}(z)=0\big)\le \exp\big( -\kappa_{in} a^2\log(n)\big).
\ee
Indeed, since the comb is a tree, covering 
$\partial\D(n-L)\cap\{z:\ |X_z|<n-AL\}$ by the DLA
cluster implies that $\D(n-AL)$ is entirely covered.
Henceforth, consider $z\in \partial \D(n-L)$ with $|X_z|\le n-AL$. 
For $y\in \D(n-AL)$, define
$h_z(y)=P_y\big(S(\EE(n-L))=z\big)$, and set
\be{def-hz}
\begin{split}
\mu(z)=&E[M_{\D(n-L)}(d(n)\ind_0;z)]-E[M_{\D(n-L)}(\D(n-L);z)]\\
=&\big(|\D(n)|-|\D(n-L)|\big) h_z(0)+\MV\big(h_z,\D(n-L)\big).
\end{split}
\ee
Using \reff{green-boundary} of
Lemma \reff{lem-green} and Lemma~\reff{lem-MV}, we obtain
when $n$ is large enough, and for $c_1,c_2>0$
independent of $n$ and $z$,
\be{mu-lower}
c_2L\big(n-L-|X_z|\big)\le \mu(z)\le c_1 L\big(n-L-|X_z|\big).
\ee
Note also that by \reff{G2-main} of Corollary \ref{cor-G2},
there are $\kappa_1,\kappa_2$ independent of $n$ and $z$
such that
\[
\kappa_2  \big(n-L-|X_z| \big)^2\le
\sum_{y\in \D(n-L)} h^2_{z}(y)\le \kappa_1 \big(n-L-|X_z| \big)^2.
\]
In order to use Lemma~\ref{lem-D}, we form the following partition
of $\D(n-L)$:
\be{def-partition}
\BB=\{y\in \D(n-L):\ X_z=X_y,\ Y_z\times Y_y\ge 0\},\quad\text{and}
\quad \A=\D(n-L)\bs \BB.
\ee
We need to show that there is $\kappa<1$, independent of $n$
such that for $y\in \A$, we have $h_z(y)<1-1/\kappa$. We choose
here $\kappa=1/2$ for simplicity, and note that for $y\in\A$,
any path joining $y$ and $z$ crosses $(X_z,0)$ so that $h_z(y)
\le h_z(X_z,0)$. Note also that
\be{ineq-A}
h_z(X_z,0)=\frac{G_{\D(n-L)}\big((X_z,0);z\big)}{2}=
\frac{G_{\D(n-L)}(0;z)}{2P_0(H(X_z,0)<\EE(n-L))}.
\ee
Using now Lemma~\ref{lem-green} and 
Proposition~\ref{prop-hitting}, we obtain readily that
$h_z(X_z,0)$ can be made smaller than 
$1/2$ as soon as $n$ is large enough.

For $\delta$ to be fixed later, we choose $A=1+\frac{1}{\delta}$ so that
$\delta\big(n-L-|X_z|\big)\ge L$, and we choose $\lambda$ as follows.
\be{choice-lambda}
\lambda=\frac{\mu(z)}{2\big(16|\BB|+\frac{1}{2}
\sum_{y\in \D(n-L)} h^2_{z}(y)\big)}\le
\frac{c_1}{2(16/3+\kappa_2)}\frac{L(n-L-|X_z|)}{(n-L-|X_z|)^2}\le
\frac{c_1}{2(16/3+\kappa_2)}\delta .
\ee
Since, we need $\lambda<\log(2)$ in Lemma~\ref{lem-D}, 
the condition on $\delta$ is that 
\be{lambda-cond}
\delta<2\log(2)\frac{16/3+\kappa_2}{c_1}.
\ee
We use Lemma~\ref{lem-D}, with $\xi=0$, $\kappa=1/2$ and 
$\lambda\le \log(2)$ to have
\[
P\big(W_{\D(n-L)}(z)=0\big)\le\exp\Big(-\lambda \mu(z)+\frac{\lambda^2}{2}
\big(\mu(z)+16|\BB|+\frac{1}{2}\sum_{y\in \D(n-L)}h^2_{z}(y)\big) \Big).
\]
Note that since $\lambda<1$, we have $\lambda \mu(z)-\lambda^2\mu(z)/2>
\lambda\mu(z)/2$, and the choice of $\lambda$ in \reff{choice-lambda}
yields \reff{inner-one} with $\kappa_{in}$ given by
$c_2/(8A^2(16/3+\kappa_1))$.

\epr

\section{Outer Fluctuations.}\label{sec-outer}
We estimate the probability that the largest finger reaches
$\partial \D(n+A\sqrt{\log(n)})$ for some large $A$. 
The analysis is distinct whether the finger protrudes in the {\it tip} 
of $\partial \D(n+A\sqrt{\log(n)})$, that is the region
\be{def-tip}
\T=\acc{z:\ |X_z|>n+\frac{A}{2}L}\cap \partial \D(n+AL),
\quad(\text{here $L:=[\sqrt{\log(n)}]$})
\ee
or in the complement of $\T$ in $\partial \D(n+AL)$
called the {\it bulk}, and made of the edges of long teeth.
Indeed, the geometry of the graph is different on the tip, 
and on the edge of a long tooth. The goal is to show that 
the appearance of a long finger
implies that a {\it narrow region} is crossed by {\it many explorers}.
More precisely, when a finger reaches site $z$ of the bulk, 
that is through a long tooth, 
it imposes that many explorers settle in the tooth:
if $Y_z>2k$ and we set $\tilde z=\A^k(z)$, then in order to cover $z$,
we need that $k$ explorers cross $\tilde z$. Moreover, at least
half of these explorers if they were random walks starting on $z$
would very likely exit $\D(n+AL)$ from $z$.
On the other hand, when a finger reaches
a site $z$ of the tip, say on $z=(n+AL,0)$, this imposes that site
$(n,0)$ is crossed by $AL$ explorers, but these explorers if they
where random walks would have many ways to exit $\D(n+AL)$.
We call $\BB(z)$ the event that
$A(n)\not\subset \D(n+AL)$, and $z$
is the first site of $\partial \D(n+AL)$
to be covered by the aggregate. Note that
\be{decomposition-A}
\acc{A(n)\not\subset \D(n+AL)}=\bigcup_{z\in \partial \D(n+AL)} \BB(z).
\ee
\subsection{The bulk region}
We assume for $z\in \partial \D(n+AL)$ with $|X_z|<n+AL/2$ that
$\BB(z)$ holds.
This implies that explorers fill the tooth of $z$
without escaping $\D(n+AL)$. For simplicity, 
we denote $\bar D=\D(n+AL)$. Let $\tilde z$ be the site of
$\partial \D(n+3AL/4)$ on the same tooth as $z$. Necessarily,
the number of explorers crossing $\tilde z$ before escaping $\bar D$
is larger than the length of the tooth to be covered
\be{def-crossing}
\xi(\tilde z)=\frac{1}{3}\big(n+AL-|X_z|\big)^2-
\frac{1}{3}\big(n+\frac{3AL}{4}-|X_z|\big)^2\ge \frac{AL}{6}
\big(n+\frac{3AL}{4}-|X_z|\big).
\ee
To exploit this information,
we introduce an auxiliary process which proved useful in studying DLA 
\cite{AG1,AG2}. The {\it flashing} process is a cluster growth,
were explorers, called $*$-explorers settle less often than in DLA.
We now build a flashing process adapted to our purpose.
\begin{itemize}
\item Inside $\D(n)$, $*$-explorers are just explorers.
\item When a $*$-explorer exits $\D(n)$, it cannot settle in $\D(n)$
anymore.
\item A $*$-explorer does not settle in $\bar D\bs \D(n)$.
\item Outside $\bar D$, $*$-explorers behave like explorers.
\end{itemize}
We call $A^*(n)$ the cluster made by $d(n)$ $*$-explorers sent
at the origin. Note that by construction, the cluster made
by $d(n)$ explorers before they exit $\D(n)$, denoted $A^*_{\D(N)}(n)$,
is equal to $A_{\D(N)}(n)$.

The key fact, established in \cite{AG1}, is that this growth
can be coupled with the internal DLA cluster in such a way
that there are times $\{T_i,T_i^*,i=1,\dots,d(n)\}$ 
with $T_i\le T_i^*$ and such that for $d(n)$ independent
random walks $S_1,\dots,S_{d(n)}$
\be{def-coupling}
A(n)=\{S_i(T_i),i=1,\dots,d(n)\},\quad\text{and}\quad
A^*(n)=\{S_i(T_i^*),i=1,\dots,d(n)\}.
\ee
Thus, under the coupling of \cite{AG1}, if an explorer happens to visit
$\tilde z$ before escaping $\bar D$, then this will be the
case for the associated $*$-explorer. We add an index $*$ to
denote objects linked with $*$-explorers. For instance, we denote
by $W^*_\Lambda(\eta,z)$ the number of $*$-explorers which
cross $z$ before escaping $\Lambda$, and we drop the $\eta$
dependence when $\eta=d(n)\ind_0$. As a consequence
of the coupling, we have
\be{coupling-1}
\BB(z)\subset \acc{W^*_{\bar D}(\tilde z)>\xi(\tilde z)}.
\ee
Let us now estimate how many
$*$-explorers exit $\bar D$ most likely from $z$.

Note first that when $|X_z|< n+AL/2$, then
\be{coupling-7}
\frac{Y_{\tilde z}}{Y_z}=
\frac{\frac{1}{3}(n+3AL/4-|X_z|)^2}{\frac{1}{3}(n+AL-|X_z|)^2}\ge
\frac{1}{4},\quad \text{so $P_{\tilde z}($ walk hits $z$ before 
$(X_z,0))$}>\frac{1}{4}.
\ee
$W^*_{\bar D}(z)$ represents the number of $*$-explorers which
exit $\bar D$ from $z$, out of $W^*_{\bar D}(\tilde z)$ $*$-explorers
at $\tilde z$. Inside $\bar D$, the $*$-explorers are just simple
random walks, and by \reff{coupling-7}, we have that 
$E[W^*_{\bar D}(z)|\ W^*_{\bar D}(\tilde z)]\ge \frac{1}{4}
W^*_{\bar D}(\tilde z)$. Thus, by Chernoff's inequality 
\be{outer-LD}
P\big(W^*_{\bar D}(z)<\frac{1}{8} W^*_{\bar D}(\tilde z)
\big| \ W^*_{\bar D}(\tilde z)\ >\xi(\tilde z)\big)\le
\exp\big(- \frac{1}{4} \frac{\xi(\tilde z)}{8}\big).
\ee
Second, note that since $z$ is a bulk site $n+3AL/4-|X_z|>AL/4$, and
on the event $\BB(z)$ we have that
$W_{\bar D}(\tilde z)> \xi(\tilde z)$, which in turn implies
$W^*_{\bar D}(\tilde z)> \xi(\tilde z)> (AL)^2/24$. Thus, we have
\be{outer-main}
P\big(W^*_{\bar D}(z)<\frac{1}{8} W^*_{\bar D}(\tilde z)
\big|\  \BB(z)\big)\le
\exp\big(- \frac{1}{24\times 32}(AL)^2\big).
\ee
Thus,
\be{coupling-2}
\begin{split}
P\big(\BB(z)\big)\le & P\big(\BB(z),
W^*_{\bar D}(z)<\frac{1}{8} W^*_{\bar D}(\tilde z)\big)+
P\big(\BB(z),
W^*_{\bar D}(z)\ge \frac{1}{8} W^*_{\bar D}(\tilde z)\big)\\
\le &P\big(W^*_{\bar D}(z)<\frac{1}{8} W^*_{\bar D}(\tilde z) 
\big| \ \BB(z) \big)
+P\big(\D(n-a_0L)\not\subset A^*_{\D(n)}(n)\big)\\
&+P\big(W^*_{\bar D}(z)\ge
\frac{\xi(\tilde z)}{8},\ \D(n-a_0L)\subset A^*_{\D(n)}(n)\big).
\end{split}
\ee
The probability $\D(n-a_0L)\not\subset A^*_{\D(n)}(n)$ is actually
estimated in Proposition~\ref{prop-inner} since 
$A^*_{\D(n)}(n)=A_{\D(n)}(n)$.

We explain now why $\{W^*_{\bar D}(z)\ge \xi\}$ 
is very unlikely, where we set for simplicity $\xi=\xi(\tilde z)/8$. 
Since our inner error estimate is also valid for $*$-explorers,
we have the equality in law
\be{upper-1}
W^*_{\bar D}(z)+M_{\bar D}\big(A^*_{\bar D}(n),z\big)
=M_{\bar D}(d(n)\ind_0,z).
\ee
This implies that
\be{upper-2}
\ind_{\D(n-a_0L)\subset A^*_{\bar D}(n)}\Big(W^*_{\bar D}(z)+
M_{\bar D}\big(\D(n-a_0L),z\big)\Big)\le M_{\bar D}(d(n)\ind_0,z).
\ee
Now, \reff{upper-2} allows us to estimate the probability that
$W^*_{\bar D}(z)$ is large, through Lemma 2.5 of \cite{AG2}:
for $0<\lambda<\log(2)$, and $\xi> 
\mu^*(z):=E[M_{\bar D}(d(n)\ind_0,z)]-E[M_{\bar D}\big(\D(n-a_0L),z\big)]$,
\be{key-outer-estimate}
P\big(W^*_{\bar D}(z)>\xi,\ \D(n-a_0L)\subset A^*_{\bar D}(n)\big)\le
\exp\Big(-\lambda(\xi-\mu^*(z))+\lambda^2\big(\mu^*(z)+4\sum_{y
\in \bar D}h_z^2(y)\big)\Big),
\ee
where $h_z(y)$ is the probability of exiting $\bar D$ from $z$, when
a random walk starts on $y$. Note that the function $y\mapsto h_z(y)$
is harmonic on $\bar D$, and that since $(A-a_0)L\ge R_{HS}$
Remark~\ref{rem-MV} applies. There is a constant $c^*$ such that
(recall that $z$ is in the bulk)
\be{upper-5}
\begin{split}
\mu^*(z)=&\big( |\D(n)|-|\D(n-a_0L)|\big) h_z(0)+\MV(h_z,\D(n-a_0L))\\
\le & c^* a_0L (n+AL-|X_z|)
\end{split}
\ee
We choose $A$ large enough, after $a_0$ is fixed, 
so that 
\be{upper-7}
\mu^*(z)\le c^* a_0L(n+AL-|X_z|)\ll
\frac{1}{4}\xi\le \frac{1}{32} \frac{AL}{4}(n+AL-|X_z|).
\ee
Also, by \reff{G2-main} of Corollary \ref{cor-G2}
we have
\be{outer-correlation}
\sum_{y\in \bar D} h_z^2(y)\le \kappa_2(n+AL-|X_z|)^2.
\ee
In the bulk, the following choice of $\lambda$ with the estimate
\reff{outer-correlation} yields
\be{def-lambda}
\lambda=\frac{\mu^*(z)}{16 \sum_{y\in \bar D} h_z^2(y)}\le
\frac{a_0 L}{16\kappa_2 (n+AL-|X_z|)}\le \frac{a_0}{8\kappa_2 A}.
\ee
One chooses $A$ large enough so that $\lambda<\log(2)$. Note that
since $\lambda<1$, our choice of $A$ in \reff{upper-7} is such that
$\lambda(\xi-\mu^*(z))-\lambda^2\mu\ge \lambda \xi/2$.
Using \reff{key-outer-estimate} with the choice of $\lambda$
in \reff{def-lambda}, and after simple algebra, 
one obtains for some constant $\kappa$ 
\be{outer-last}
P\big(W^*_{\bar D}(z)>\xi,\ \D(n-a_0L)\subset A^*_{\bar D}(n)\big)\le
\exp\big(-\kappa a_0 A L^2\big).
\ee
Finally, from the inner error, we know that for $a_0$ large enough
most likely $\D(n-a_0L)\subset A_{\D(n)}(n)=A^*_{\D(n)}(n)\subset
A^*_{\bar D}(n)$. Combining the estimates for the three terms
on the right hand side of \reff{coupling-2}, we obtain
\[
P\big(\BB(z)\big)\le \exp\big(-\frac{1}{24\times 32}A^2\log(n)\big)+
\exp\big(-\kappa a_0A\log(n)\big)+\exp\big(-\kappa_{in}a^2_0\log(n)\big).
\]
We can choose $a_0$, and then $A$ so that $P(\BB(z))$ is smaller
than any negative power of $n$.
\subsection{The Tip.}
Let us describe the additional idea needed to deal with the tip.
A constant $A$ large enough will be chosen later. Assume $A/4\in \N$,
and define three points
\[
B=(n+\frac{A}{2}L,0),\quad \tilde B=(n+\frac{A}{4}L,0),\quad
\text{and}\quad C\in \partial \D(n+AL)\text{ with } X_C=n+
\frac{A}{4}L.
\]
Assume in this section that a site of the tip is covered.
This implies that $B$ or $-B$ is in $A(n)$. Assume 
for instance that $B\in A(n)$.
The internal DLA covering mechanism would say that 
$\tilde B$ is necessarily covered by $\frac{A}{4}L$ explorers.
However, too small a fraction, of the order of $1/L$, 
of these explorers 
would exit $\D(\tilde B, \frac{3}{4}AL)$ from site $C$.
We first need to show that of the order
of $(AL)^3$ explorers cross $\tilde B$, and secondly that
it is very unlikely that of the order of $(AL)^2$ exit 
$\D(\tilde B, \frac{3}{4}AL)$ from site $C$.

As in the previous section, we need to consider 
here the same $*$-explorers.
An important property is the fact that the
aggregate's law is independent of the order of the explorers
we launch, or more generally, of stopping explorers
in some region letting other explorers cover space before the stopped
ones are eventually launched. Thus, we will realize the aggregate by sending
two waves of exploration. We stop $*$-explorers on $\tilde B\cup
\partial D(n+AL)$, and call $\zeta$ the configuration of stopped
$*$-explorers, that is $\zeta:\{\tilde B\}\cup \partial \D(n+AL)\to \N$.

The event that $B$ is covered, and
$\zeta(\tilde B)$ is less than $\beta (AL/4)^3$, is very unlikely by
Corollary~\ref{cor-cross}. 
Henceforth, assume that $\zeta(\tilde B)>\beta (AL/4)^3$,
where we recall that $\beta$ is a constant independent of $A,n$.
Assume that we launch the stopped $*$-explorers and stop them
on $\partial \D(n+AL)$. It is very unlikely that 
less than $\kappa (AL)^2$ $*$-explorers 
exit $\partial \D(n+AL)$ from $C$ for some positive constant
$\kappa$. Indeed, let us call the latter number $W^*_{\bar D}
\big(\zeta(\tilde B)\ind_{\tilde B};C\big)$. Call for simplicity
$\Lambda=\D(\tilde B, \frac{3}{4}AL)$, and define
$M_{\Lambda}\big(\zeta(\tilde B)\ind_{\tilde B};C\big)$
the number of random walks which exit $\Lambda$
from $C$. Note the following obvious fact 
\[
\Lambda\subset \bar D\Longleftrightarrow
W^*_{\bar D} \big(\zeta(\tilde B)\ind_{\tilde B};C\big)\ge
W^*_{\Lambda} \big(\zeta(\tilde B)\ind_{\tilde B};C\big),
\]
Also, the abelian property we mentioned says that (with equality in law)
\be{cor-abelian}
W^*_{\bar D}(C)=
W_{\bar D}^* \big(\zeta(\tilde B)\ind_{\tilde B};C\big),
\ee
Thus, to estimate the probability that $W^*_{\bar D}(C)$ be small
it is enough to estimate the probability that $W^*_{\Lambda}
(\zeta(\tilde B)\ind_{\tilde B};C)$ be small.
Note that $*$-explorers when starting in $\tilde B$ and staying
in $\Lambda$ have the same law as simple random walks, so that 
\[
W^*_{\Lambda}\big(\zeta(\tilde B)\ind_{\tilde B};C\big)=
M_{\Lambda}\big(\zeta(\tilde B)\ind_{\tilde B};C\big).
\]
$M_{\Lambda}$ is a sum of independent Bernoulli, and from
\reff{green-boundary}, there is a constant $\kappa>0$
\[
E\cro{M_{\Lambda}\big(\zeta(\tilde B)\ind_{\tilde B};C\big)
\big| \zeta(\tilde B)\ge \beta(\frac{AL}{4})^3}\ge
\beta(\frac{AL}{4})^3\frac{2}{3AL/4}\ge 2\kappa (AL)^2,
\]
and therefore, using Chernoff's inequality
on the event $\{\zeta(\tilde B)\ge \beta(\frac{AL}{4})^3\}$
\be{estimates-M}
\begin{split}
P\big( W^*_{\Lambda}\big(\zeta(\tilde B)\ind_{\tilde B};C\big)
\le \kappa (AL)^2 \big)=&
P\big(M_{\Lambda}\big(\zeta(\tilde B)\ind_{\tilde B};C\big)
\le \kappa (AL)^2 \big)\\
\le& \exp\big(-\kappa (AL)^2\big).
\end{split}
\ee
We deal now with the event $\{W^*_{\bar D}(C)>\kappa (AL)^2\}$.
Note that defining
\[
\mu^*_{\bar D}(C):=E[M_{\bar D}(d(n)\ind_0,C)]-
E[M_{\bar D}(\D(n-a_0L),C)],
\]
we have using our harmonic measure estimate, for some constant $c^*$
\be{tip-1}
\mu^*_{\bar D}(C)\le c^* a_0AL^2,\quad
\text{and}\quad \sum_{y\in \D(n-a_0L)}h^2_C(y)\le
\sum_{y\in \bar D} h^2_C(y)\le\kappa_2(\frac{3AL}{4})^2.
\ee
Thus, for $\lambda<\log(2)$, 
\be{tip-2}
P\big(W^*_{\bar D}(C)>\kappa (AL)^2,\ \D(n-a_0L)\subset
A^*_{\bar D}(n)\big)\le 
e^{-\lambda(\kappa (AL)^2-c^*a_0AL^2)+\lambda^2
(a_0AL^2+\kappa_2 (AL)^2)}.
\ee
We need now to choose $A$ so large that $\kappa A\ge 2 \kappa_2 c^*a_0$, and
$\lambda=\min(\log(2),\ \kappa/(2\kappa_2))$, which gives
finally that
\be{tip-main}
P(\text{the tip is covered})\le \exp\big(-c A^2L^2\big).
\ee
\vskip 2,cm
\centerline{\large \bf APPENDIX}
\vskip ,5cm
\appendix
\section{Proof of Lemma~\ref{lem-green}}
Our goal in this section is to estimate precisely the restricted
Green's function
$z\mapsto G_{\D(\rho)}(0;z)$ for any positive
real $\rho$, and $z\in \D(\rho)$.
We use that the latter function is discrete harmonic on
$\D(\rho)\bs\{0\}$, vanishes on the (discrete)
boundary of $\D(\rho)$, and satisfies $\Delta G_{\D(\rho)}(0,\cdot)|_0=-1$.

We first find an explicit function, denoted $h:\Z\times \R$,
discrete harmonic on the $x$-axis,
real harmonic on each tooth of $\D(\rho)\bs\{0\}$,
vanishing on
\[
\Sigma(\rho)=\acc{(x,y)\in \Z\times \R:\
|x|\le \rho ,\ |y|=\frac{1}{3}\pare{\rho-|x|}^2},
\quad\text{and}\quad \Delta h(0)=-1.
\]
Since $h$ is linear on each tooth of $\D(\rho)$, and can readily
be extended to $\D(\rho)\cup\partial \D(\rho)$ with non-positive
values on $\partial \D(\rho)$, the maximum principle yields
\be{estimate-1}
\forall z\in \D(\rho)\qquad G_{\D(\rho)}(0;z)\ge h(z).
\ee
Similarly we build $h^+$, positive and harmonic on a larger domain
$\Sigma^+(\rho)\bs\{0\}$ with
\[
\Sigma^+(\rho)=\acc{(x,y)\in \Z\times \R:\
|x|\le [\rho]+1 ,\ |y|=\frac{1}{3}\pare{[\rho]+1-|x|}^2},
\quad\text{and}\quad \Delta h^+(0)=-1.
\]
We will have that $h^+-G_{\D(\rho)}$
is harmonic on $\D(\rho)$ and nonnegative on $\partial \D(\rho)$.
Again, by the maximum principle
\be{estimate-2}
\forall z\in \D(\rho)\qquad h^+(z)\ge G_{\D(\rho)}(0;z).
\ee
The explicit expression of $h$ and $h^+$, and
estimates \reff{estimate-1} and \reff{estimate-2} are
the desired results of this section.

\paragraph{Construction of $h$.}
By the symmetries
of $\D(\rho)$, $h$ is even in the $x$ and $y$ coordinates. Thus,
we restrict the construction for $x\in [-\rho,0]$.
Also, it is convenient to shift $\D(\rho)$ by $\rho$
units along the $x$-axis, so that $(-\rho,0)$ becomes the origin,
and $(0,0)$ becomes $(0,\rho)$.

On each arm of the comb $h$ is linear, and reads for $z=(x,y)$,
\be{step-1}
h(z)= a(x)y+b(x).
\ee
We set $Y(x)=\frac{1}{3}x^2$, and we impose
\be{step-2}
0=a(x)Y(x)+b(x),\quad\text{and}\quad b(0)=0.
\ee
We solve a set of equations: for $x_k=\rho-[\rho]+k$ with
$k$ integer in $\{1,\dots,[\rho]-1\}$,
\be{step-3}
4h(x_k,0)=h(x_k,1)+h(x_k,-1)+h(x_{k+1},0)
+h(x_{k-1},0)=2h(x_k,1)+h(x_{k+1},0)+h(x_{k-1},0),
\ee
and a boundary equation
\be{step-4}
4h(\rho,0)=2h(\rho,1)+2h(\rho-1,0)+4.
\ee
When we choose $b(x)=\frac{1}{3}\alpha x^3$, \reff{step-2} implies that
$a(x)=-\alpha x$.
In terms of $a$ and $b$,
\reff{step-3} and \reff{step-4} read for $k\in \{1,\dots,[\rho]-1\}$
\be{step-5}
2b(x_k)=2a(x_k)+b(x_{k+1})+b(x_{k-1}),\quad\text{and}\quad
b(\rho)=a(\rho)+b(\rho-1)+2.
\ee
Solving \reff{step-5}, we find
\[
\alpha=\frac{2}{\rho^2+\frac{1}{3}}.
\]
Thus, we obtain a function $h:[0,\rho]\times\R$ given by
\be{def-tilde-h}
h(x,y)=\frac{2x}{\rho^2+\frac{1}{3}}\big(\frac{x^2}{3}-y\big).
\ee
\paragraph{Construction of $h^+$.}
Here, the domain $\D([\rho]+1)$ is shifted by $[\rho]+1$
units along the $x$-axis.
We build a function $h^+(x,y)=a^+(x)y+b^+(x)$,
with $h^+(0,0)=0$, and $h^+(x,Y^+(x))=0$ for
\[
Y^+(x)=\frac{1}{3}x^2+1,\quad\text{and}\quad
b^+(x)=\alpha x Y^+(x).
\]
This implies that $a^+(x)=-\alpha x$. Now, $a^+,b^+$
solve \reff{step-3} for $x$ an integer
from 1 to $[\rho]+1$. Also \reff{step-4} holds with $[\rho]+1$
instead of $\rho$. This yields
\[
\alpha=\frac{2}{([\rho]+1)^2+\frac{1}{3}+1}.
\]
We obtain a function $h^+:[0,[\rho]+1]\times\R$ given by
\be{def-tilde-h+}
h^+(x,y)=\frac{2x}{([\rho]+1)^2+\frac{1}{3}+1}\big(\frac{x^2}{3}+1-y\big).
\ee
\paragraph{Estimate on Green's function.}
We go back now to the usual coordinate system to obtain
\reff{green-asymp} with $h$ and $h^+$ given in \reff{def-h} and
\reff{def-h+}.
One more relation is useful to obtain nondegenerate estimates
when taking $z$ close to the boundary of $\D(\rho)$. Recall that
for $z\not= 0$,
\be{A-2}
\begin{split}
G_{\D(\rho)}(0;z)=& E_0\big[\sum_{k=1}^{\EE(\rho)-1}\ind_{S(k)=z}\big]
= E_0\big[\sum_{k=1}^{\EE(\rho)-1}\ind_{S(k-1)=\A(z),\ S(k)=z}\big]\\
=&G_{\D(\rho)}(0;\A(z))\times p\big(\A(z);z\big).
\end{split}
\ee
Thus, if $z=(x,0)\in \D(\rho)$, and $x+1>\rho$, note that
$\A(z)=(x-1,0)$, and $p\big(\A(z);z\big)=1/4$ so that
\be{A-3}
G_{\D(\rho)}(0;z)=\frac{1}{4}G_{\D(\rho)}(0;\A(z))\quad\text{and}\quad
G_{\D(\rho)}(0;z)\ge
\frac{h\big((x-1,0)\big)}{4}.
\ee
\section{Proof of Proposition~\ref{prop-hitting}}
Proposition~\ref{prop-hitting} uses the following Lemma
which we prove at the end of the section.
\bl{lem-1step}
For any $\rho>0$,
\be{ineq-1step}
P_0\big(H(1,0)<\EE(\rho)\big) \ge \big(\frac{\rho}{\rho+1}\big)^3.
\ee
\el

We assume that the integer $y$ satisfies $0<y<\n-1$,
and denote for $k=0,\dots,y$
\be{def-a}
u(k)=P_{(k,0)}\big( H(y,0)<\EE(\n)\big).
\ee
and $L_\n(k)$ denotes the height of the tooth
of $\D(\n)$ at site $(k,0)$, i.e.
$L_\n(k)=[\frac{1}{3}(\n-k)^2]+1$.
If we condition the event
$\{H(y,0)<\EE(\n)\}$ on the first step of the random walk,
then we obtain for $k=1,\dots,y-1$
\be{acc-1}
u(k)=\frac{u(k+1)+u(k-1)}{4}+\frac{1}{2}\big(1-
\frac{1}{L_\n(k)}\big) u(k),\quad\text{and}\quad
u(y)=1.
\ee
We rewrite \reff{acc-1} as
\be{acc-2}
u(k)-\frac{u(k-1)}{\alpha(k-1)}=\frac{1}{\alpha(k-1)}
\big(u(k+1)-\frac{u(k)}{\alpha(k)}\big),
\ee
and the $\{\alpha(k),\ k=0,\dots,y-1\}$ is a sequence
obtained inductively with the constraint that
\be{def-alpha}
\forall k=1,\dots,y-1\qquad
\alpha(k-1)+\frac{1}{\alpha(k)}=2+\frac{2}{L_\n(k)}.
\ee
As we iterate \reff{def-alpha}, from $k=1$ to $k=y-1$, we find
\be{acc-3}
u(1)-\frac{u(0)}{\alpha(0)}=\frac{1}{\alpha(0)}\times
\dots\times \frac{1}{\alpha(y-2)}\big( u(y)-
\frac{u(y-1)}{\alpha(y-1)}\big).
\ee
Now, assume we have the following three relations
\be{require-3}
(i)\quad u(0)\ge u(1) \big(\frac{\n}{\n+1}\big)^3\qquad
(ii)\quad 1-\frac{u(y-1)}{\alpha(y-1)}\ge \frac{1}{\n-y-1},
\ee
and
\[
(iii)\quad 1\le \alpha(k)\le 1+\frac{3}{\n-k-2}\quad \forall k<y.
\]
Using \reff{acc-3} and \reff{require-3}, we obtain
for $\n$ large enough and for a positive constant $\kappa$
\be{acc-4}
\big((\frac{\n+1}{\n})^3-\frac{\n-2}{\n+1}\big) u(0)\ge
\Big(\prod_{k=0}^{y-2} \frac{\n-k-2}{\n-k+1}\Big)\times
\Big(1-\frac{u(y-1)}{\alpha(y-1)}\Big)
\ee
Now, some simple algebra yields
\be{A-5}
\big((\frac{\n+1}{\n})^3-\frac{\n-2}{\n+1}\big)\le
\frac{6}{\n}\big(1+\frac{1}{2\n}+\frac{1}{6\n^2})\big),
\ee
and
\be{A-6}
\prod_{k=0}^{y-2} \frac{\n-k-2}{\n-k+1}\ge \frac{(\n-y)^3}{\n^3}.
\ee
We deduce from \reff{acc-4}, \reff{A-5} and \reff{A-6} that
for some constant $\kappa$
\be{A-7}
u(0)\ge \kappa \big(\frac{\n-y}{\n}\big)^2.
\ee
We are left with showing the estimates of \reff{require-3}.
Note that (i) is Lemma~\ref{lem-1step}.

We show (ii). When we start on $(y-1,0)$, one way to
to escape $\D(\n)$ before reaching $(y,0)$ is to go up
on one tooth and hit the boundary of $\D(\n)$ before touching $(y-1,0)$.
Thus, when $y<\n-1$
\be{acc-5}
1-u(y-1)=P_{(y-1,0)}(\EE(\n)<H(y,0))\ge
\frac{1}{2L_\n(y-1)}\ge \frac{3/2}{(\n-y+1)^2+1}\ge \frac{1}{(\n-y+1)^2}.
\ee
This is equivalent to
\be{acc6}
1-\frac{u(y-1)}{1+\frac{1}{\n-y+1}}\ge \frac{1}{\n-y+1}.
\ee
Now, to produce a sequence satisfying \reff{def-alpha}, we
choose $\alpha(y-1)$ as follows, and build $\alpha(k)$ by
a backward induction:
\be{def-alpha0}
\alpha(y-1)=1+\frac{3}{\n-y+1}.
\ee
This implies, using \reff{acc6}, that
\be{acc-7}
1-\frac{u(y-1)}{\alpha(y-1)}\ge 1-\frac{u(y-1)}{1+\frac{1}{\n-y+1}}
\ge \frac{1}{\n-y+1}.
\ee

We now show that (iii) is compatible with our choice
\reff{def-alpha0}. We do it by backward induction.
First, it is obvious that $\alpha(k)>1$ implies
that $\alpha(k-1)>1$. We assume now that $\alpha(k-1)> 1
+3/(\n-k-1)$, and show that $\alpha(k)> 1+3/(\n-k-2)$. This
in combination with \reff{def-alpha0} yields (iii).
In view of \reff{def-alpha} this is equivalent to checking that
\be{acc-8}
\begin{split}
& 1\ge \big( 1+\frac{2}{L_\n(k)}-\frac{3}{\n-k-1}\big)
\big( 1+\frac{3}{\n-k-2}\big)\\
\Longleftrightarrow &
1\ge 1-\frac{9}{(\n-k-1)(\n-k-2)}+\big(\frac{3}{\n-k-2}-\frac{3}{\n-k-1}\big)+
\frac{2}{L_\n(k)}\big( 1+\frac{3}{\n-k-2}\big)\\
\Longleftrightarrow &
L_\n(k)\ge \frac{(\n-k-1)(\n-k+1)}{3}=\frac{1}{3}\big((\n-k)^2-1\big).
\end{split}
\ee
The last inequality of \reff{acc-8} is true since
$L_\n(k) \ge \frac{1}{3}(\n-k)^2$.

\noindent{\bf Proof of Lemma~\ref{lem-1step}}

Calling $A=(1,0)$, we establish first,
\be{step-9}
P_0(H(A)<\EE(\n))=\big( 1+\frac{1}{L(\n)}+\frac{2}{G_{\D(\n)}(0;0)}
\big)^{-1}.
\ee
For simplicity, we
name $B=(0,1)$, $C=(-1,0)$ and $D=(0,-1)$. Then,
\be{step-6}
\begin{split}
P_0(H(A)<\EE(\n))&=P_0(S_1=A)+\sum_{z\in\{B,C,D\}} P_0(S_1=z)P_z
\big(H(0)<\EE(\n)\big)P_0(H(A)<\EE(\n))\\
&= \frac{1}{4}+P_0(H(A)<\EE(\n))\big(\frac{P_B(H(0)<\EE(\n))}{2}+
\frac{ P_C(H(0)<\EE(\n))}{4}\big).
\end{split}
\ee
A classical gambler's ruin estimate yields
\be{step-7}
P_B(H(0)<\EE(\n))=\frac{L_\n(0)-1}{L_\n(0)}.
\ee
Also, by decomposing over the first step, we have for $H(0)^+$ the
return time to 0,
\be{step-8}
\begin{split}
1-\frac{1}{G_{\D(\n)}(0;0)}=& P_0(H(0)^+<\EE(\n)) \\
=& \frac{1}{2} P_A(H(0)<\EE(\n))+\frac{1}{2} P_B(H(0)<\EE(\n)).
\end{split}
\ee
Thus, using \reff{step-7} and \reff{step-8} in \reff{step-6},
we obtain \reff{step-9}.
Now, by \reff{green-asymp}, we have
\be{new-5}
\frac{2}{G_{\D(\n)}(0;0)}\le \frac{3}{\n}+\frac{1}{\n^3}.
\ee
Recalling that $L_\n(0)\ge \n^2/3$, we obtain the desired relation.
\be{main-P}
P_0(H(A)<\EE(\n))\ge \big(1+\frac{3}{\n}+
\frac{3}{\n^2}+\frac{1}{\n^3}\big)^{-1}=\Big(\frac{\n}{\n+1}\Big)^3.
\ee
\section{Proof of Green's function estimates.}
\subsection{Proof of Proposition~\ref{prop-green}}
We prove (i). By symmetries of $\D(\n)$, we can consider
$0<X_w<X_z$ or $X_w=X_z$ and $Y_w,Y_z\ge 0$.
Note that the path joining $w$ and $z$ crosses $(X_w,0)$, as well
as the path joining 0 and $z$. Thus,
\be{case-1}
\begin{split}
G_{\D(\n)}(w;z)=&P_w\big(H(X_w,0)<\EE(\n)\big)\times
G_{\D(\n)}\big((X_w,0);z\big),\qquad\text{and}\\
G_{\D(\n)}(0;z)=&P_0\big(H(X_w,0)<\EE(\n)\big)\times G_{\D(\n)}
\big((X_w,0);z\big).
\end{split}
\ee
This implies that
\be{case-2}
G_{\D(\n)}(w;z)=P_w\big(H(X_w,0)<\EE(\n)\big)\times
\frac{G_{\D(\n)}(0,z)}{P_0\big(H(X_w,0)<\EE(\n)\big)}.
\ee
Since $P_w\big(H(X_w,0)<\EE(\n)\big)=(L_\n(w)-Y_w)/L_\n(w)$, \reff{green-1}
follows from Lemma \ref{lem-green} and Proposition \ref{prop-hitting}.

Note that case (ii) follows from the previous argument
by noting first that a reversible measure for the simple random
walk assigns to a vertex its degree, and
\be{sym-1}
G_{\D(\n)}(w;z)=\frac{\deg(z)}{\deg(w)} G_{\D(\n)}(z,w)
\Longrightarrow
G_{\D(\n)}(w;z)\le 2G_{\D(\n)}(z,w).
\ee
Then, we interchange in \reff{case-2}
the role of $z$ and $w$. Note, however that
$z$ is at a distance 1 from the boundary of $\D(\n)$ while
$w$ can be anywhere in $\D(\n)$.

We prove now (iii). Note the two relations
\[
G_{\D(\n)}\big((X_w,0);z\big)=
P_{(X_w,0)}\big(H(0)<\EE(\n)\big)\times G_{\D(\n)}(0;z),
\]
and,
\be{case-3}
G_{\D(\n)}\big(0;(X_w,0)\big)\le
2 G_{\D(\n)}\big((X_w,0);0\big)=
P_{(X_w,0)}\big(H(0)<\EE(\n)\big)\times G_{\D(\n)}(0,0).
\ee
Using Lemma~\ref{lem-green} and
Proposition~\ref{prop-hitting}, \reff{case-3}
yields
\be{case-4}
G_{\D(\n)}\big((X_w,0);z\big)=\frac{G_{\D(\n)}\big(0;(X_w,0)\big)}
{G_{\D(\n)}(0;0)}\times G_{\D(\n)}(0;z)\le
\kappa \frac{(\n-X_w)^3}{3\n^2+1}\times
\frac{3\n^2+1}{2\n^3}\times \frac{(\n-X_z)}{\n^2}.
\ee
We complete \reff{case-4} with the gambler's ruin estimate to obtain
\reff{green-3}.

Finally, we deal with (iv). Consider $w,z$ with $X_w=X_z$. Then,
\be{case-6}
G_{\D(\n)}(w;z)=P_w(H(z)<\EE(\n))G_{\D(\n)}(z;z)\le G_{\D(\n)}(z;z).
\ee
On the other hand,
by decomposing over the first step (and recalling that $z\in \partial_I
\D(\n)$)
\be{case-7}
G_{\D(\n)}(z;z)=1+\frac{1}{2} G_{\D(\n)}(z-1,z)\le 1+\frac{1}{2} G_{\D(\n)}(z;z)
\Longrightarrow
G_{\D(\n)}(z;z)\le 2.
\ee
This completes \reff{green-4}.
\subsection{Proof of Corollary \ref{cor-G2}.}
It is enough to consider $z\in \partial\D^+(\n)$, and to recall
the last passage decomposition \reff{last-p}.
Introduce now the following notation. For $\Lambda\subset \Z^2$,
\be{def-Gamma}
\Gamma(\Lambda)=\sum_{w\in \Lambda} G_{\D(\n)}^2(w,\A(z)),
\ee
and partition $\D(\n)$ into four regions $\D_1,\dots,\D_4$ with
\[
\D_1=\D(\n)\cap\big( \acc{0\le x\le X_z}\cup\acc{x=X_z,\ Y_z\cdot
Y_w<0}\big),\quad \D_2=\D(\n)\cap \acc{X_w=x,\ Y_z\cdot Y_w\ge 0},
\]
$\D_3=\D(\n)\cap \acc{x>X_z}$,
and $\D_4$ the remaining part of $\D(\n)$. Using $\kappa_i$ to
denote constants, whose meaning may change from line to line,
we obtain using the estimates of Proposition~\ref{prop-green}.
We set $n=[\n]+1$, and
\be{Gamma-1}
\Gamma(\D_1)\le \kappa (\n+1-X_z)^2 \sum_{k=n-X_z}^n
\frac{1}{k^4}\sum_{i=1}^{k^2} \frac{i^2}{k^4}\le
\kappa_1 (\n+1-X_z)^2,
\ee
then,
\be{Gamma-2}
\Gamma(\D_2)\le 2 L_n(z)\le \frac{2(\n+1-X_z)^2}{3}.
\ee
Also, note the lower bound
\be{Gamma-lower}
\Gamma(\D_2)\ge \kappa \sum_{i=1}^{[(\n-X_{\A(z)})^2]}
\frac{i^2}{(\n-X_{\A(z)})^4}\ge \sous \kappa (\n-X_{\A(z)})^2.
\ee
Now,
\be{Gamma-3}
\Gamma(\D_3)\le\kappa \frac{1}{(\n-X_z)^8}\sum_{k=1}^{n-X_z} k^5
\le \kappa_3 \frac{1}{(\n+1-X_z)^2}.
\ee
Finally
\be{Gamma-4}
\Gamma(\D_4)\le \frac{(\n+1-X_z)^2}{\n^{10}}\sum_{k=1}^n k^6
\sum_{i=1}^{k^2} \frac{i^2}{k^4}\le \kappa_3 \frac{ (\n+1-X_z)^2}{\n}.
\ee
With our estimates, the dominant term in $\Gamma(\D(\rho))$
is $\Gamma(\D_2)$, and this concludes the proof.

\subsection{On sums of Bernoulli.}
Let us recall Lemma 2.3 of \cite{AG1}.
Assume that for random variables $W,M$, and $L$ we have
\be{D-1}
W+L\ge M,
\ee
and furthermore that $L$
and $M$ are sums of independent Bernoulli variables with
$L=Y_1+\dots+Y_n$. Three hypotheses played a key role
in \cite{AG1}: $(\H_0)$ $W$ is
independent of $L$,
\be{hyp-D}
(\H_1)\quad \mu:=E[M]-E[L]\ge 0,\qquad\text{and}\quad
(\H_2)\quad \text{for some $\kappa>1$} \quad
\sup_i E[Y_i]<1-\frac{1}{\kappa}.
\ee
Then, Lemma 2.3 of \cite{AG1} establishes that for $0\le \xi<\mu$,
any $\lambda\ge 0$,
\be{main-AG}
P\big(W\le \xi\big)\le \exp\big(-\lambda(\mu-\xi)+
\frac{\lambda^2}{2}(\mu+\kappa\sum_i E[Y_i]^2)\big).
\ee
In the inner estimate that we treat here, hypothesis $(\H_2)$ does not hold. Rather,
we decompose the Bernoulli variables $\{Y_i,\ 1\le i\le n\}$
into two subgroups, according to some $\kappa>1$ as follows:
\be{def-AB}
\A=\{i\le n:\ E[Y_i]<1-\frac{1}{\kappa}\},\quad\text{and}\quad
\BB=\{1,\dots,n\}\bs \A.
\ee
We show the following estimate.
\bl{lem-D}
For $\{W,M,L\}$ satisfying $\H_0$ and $\H_1$,
and any $0\le \lambda\le \log(2)$, we have
\be{main-bernoulli}
\text{for}\quad \xi\ge 0,\quad
P\big(W\le \xi\big)\le \exp\Big(-\lambda(\mu-\xi)+
\frac{\lambda^2}{2}\big(\mu+\frac{4}{\kappa^2}
|\BB|+\kappa\sum_{i\in \A} E[Y_i]^2
\big)\Big).
\ee
\el
\bpr
Using Chebychev's inequality with any $\lambda>0$, and hypothesis ($\H_0$)
\[
P\big(W\le \xi\big)\le e^{\lambda \xi}
\frac{E[e^{-\lambda W}]E[e^{-\lambda L}]}{E[e^{-\lambda L}]}\le
e^{\lambda \xi} \frac{E[e^{-\lambda M}]}{E[e^{-\lambda L}]}\le
e^{-(\mu-\xi)\lambda} \frac{E[e^{-\lambda (M-E[M])}]}
{E[e^{-\lambda(L-E[L])}]}.
\]
We have now to estimate the Laplace transform of Bernoulli variables.
The argument follows the proof of Lemma 2.3 of \cite{AG1}, with
the following trick. When $i\in \BB$, $\tilde Y_i=1-Y_i$
is again a Bernoulli variable, and
\be{trick-1}
Y_i-E[Y_i]=-\big(\tilde Y_i-E[\tilde Y_i]\big).
\ee
Now, we recall two simple inequalities used in the
proof of Lemma 2.3 of \cite{AG1}: for $0\le x\le 1$, we have
$1+x\ge \exp(x-x^2)$, whereas for $0\le x\le 1-1/\kappa$, we have
$1-x\ge \exp(-x-\kappa x^2/2)$.
Thus, using $e^\lambda\le 2$ and the notation $f(t)=e^t-1-t$, and
$g(t)=(e^t-1)^2$, we have for $i\in \BB$,
\be{trick-2}
\begin{split}
E[\exp\big(-\lambda(Y_i-E[Y_i])\big)]=&E[\exp\big(\lambda(
\tilde Y_i-E[\tilde Y_i]\big)]=e^{-E[\tilde Y_i]\lambda}
\big(1+E[\tilde Y_i](e^\lambda-1)\big)\\
\ge & \exp\big(f(\lambda)E[\tilde Y_i]-g(\lambda)E[\tilde Y_i]^2\big).
\end{split}
\ee
On the other hand, for $i\in \A$,
\be{same-A}
E[\exp\big(-\lambda(Y_i-E[Y_i])\big)]\ge
\exp\big(f(-\lambda)E[Y_i]-\frac{\kappa}{2} g(-\lambda)
E[Y_i]^2\big).
\ee
Recall now that if $[.]_+$ stands for the positive part
\[
0\le f(t)\le \frac{t^2}{2}e^{[t]_+},\quad\text{and}\quad
0\le g(t)\le t^2 e^{2[t]_+}.
\]
Finally, we have (using also in the third line
that for $\lambda\ge 0$, we have $f(\lambda)\ge f(-\lambda)$)
\be{trick-3}
\begin{split}
\frac{E[e^{-\lambda(M-E[M])}]}{E[e^{-\lambda(L-E[L])}]}\le &
\exp\big(f(-\lambda)E[M]-f(-\lambda)\sum_{i\in \A} E[Y_i]
-f(\lambda) \sum_{i\in \BB} E[\tilde Y_i]\\
&\qquad +
g(\lambda)\sum_{i\in \BB}  E[\tilde Y_i]^2+
\frac{\kappa}{2} g(-\lambda) \sum_{i\in \A} E[Y_i]^2\big)\\
\le& \exp\Big(f(-\lambda)\mu +(f(-\lambda)-f(\lambda))
\sum_{i\in \BB} E[\tilde Y_i]\\
&\qquad+ g(\lambda)\sum_{i\in \BB}  (1-E[Y_i])^2+
\frac{\kappa}{2} g(-\lambda) \sum_{i\in \A} E[Y_i]^2\Big) \\
\le & \exp\Big( \frac{\lambda^2}{2}\big( \mu +\frac{4}{\kappa^2}|\BB|+
\frac{\kappa}{2} \sum_{i\in \A} E[Y_i]^2\big)\Big).
\end{split}
\ee
\epr

\noindent{\bf Acknowledgements.}
The authors thank an anonymous referee for careful
reading, and suggestions which considerably 
improved the presentation.

\end{document}